\documentclass[12pt,reqno]{amsart}

\usepackage{amsthm,amsfonts,amssymb,euscript}
\numberwithin{equation}{section}
\begin{document}

\title[Schr\"{o}dinger maps]
{Low-regularity Schr\"{o}dinger maps}
\author{Alexandru D. Ionescu}
\address{University of Wisconsin--Madison}
\email{ionescu@math.wisc.edu}
\author{Carlos E. Kenig}
\address{University of Chicago}
\email{cek@math.uchicago.edu}
\thanks{The first author was
supported in part by an NSF grant and a Packard fellowship.
The second author was supported in part by an NSF grant.}
\begin{abstract}
We prove that the Schr\"{o}dinger map initial-value problem
\begin{equation*}
\begin{cases}
&\partial_ts=s\times\Delta_x s\,\text{ on }\,\mathbb{R}^d\times[-1,1];\\
&s(0)=s_0
\end{cases}
\end{equation*}
is locally well-posed for small data $s_0\in H^{{\sigma_0}}_Q(\mathbb{R}^d;\mathbb{S}^2)$, ${\sigma_0}>(d+1)/2$, $Q\in\mathbb{S}^2$.
\end{abstract}
\maketitle

\section{Introduction}\label{intro}

In this paper we consider the Schr\"{o}dinger map initial-value problem
\begin{equation}\label{Sch1}
\begin{cases}
&\partial_ts=s\times\Delta_x s\,\text{ on }\,\mathbb{R}^d\times[-1,1];\\
&s(0)=s_0,
\end{cases}
\end{equation}
where $d\geq 2$ and $s:\mathbb{R}^d\times[-1,1]\to\mathbb{S}^2\hookrightarrow\mathbb{R}^3$ is a smooth function. The Schr\"{o}dinger map equation has a rich geometric structure and arises naturally in a number of different ways; we refer the reader to \cite{NaStUh} or \cite{KePoStTo} for details. In this paper we prove a local well-posedness result for the initial-value problem \eqref{Sch1} for small data in low-regularity spaces. 

For $\sigma\geq 0$ let $J^\sigma$ denote the operator on $\mathcal{S}'(\mathbb{R}^d)$ defined by the Fourier multiplier $\xi\to(1+|\xi|^2)^{\sigma/2}$, and let $H^{\sigma}=H^\sigma(\mathbb{R}^d)$ denote the usual Banach spaces of complex-valued Sobolev functions on $\mathbb{R}^d$, $\|f\|_{H^\sigma}=\|J^\sigma(f)\|_{L^2}$. For $\sigma\geq 0$ and $Q=(Q_1,Q_2,Q_3)\in\mathbb{S}^2$ we define the complete metric space
\begin{equation}\label{Sch2}
H^{\sigma}_Q(\mathbb{R}^d;\mathbb{S}^2)=\{f:\mathbb{R}^d\to\mathbb{R}^3:|f(x)|\equiv 1\text{ and }f_l-Q_l\in H^\sigma\text{ for }l=1,2,3\},
\end{equation}
with the induced distance
\begin{equation}\label{Ban2}
d^\sigma_Q(f,g)=\big[\sum_{l=1}^3\|f_l-g_l\|_{H^\sigma}^2\big]^{1/2}.
\end{equation}
For $Q\in\mathbb{S}^2$ we define the complete metric space
\begin{equation*}
H^{\infty}_Q(\mathbb{R}^d;\mathbb{S}^2)=\bigcap_{\sigma\geq 0}H^{\sigma}_Q(\mathbb{R}^d;\mathbb{S}^2)\text{ with the induced metric.}
\end{equation*}

For $Q\in \mathbb{S}^2$ let $f_Q(x)\equiv Q$, $f_Q\in H^{\infty}_Q(\mathbb{R}^d;\mathbb{S}^2)$. For any metric space $X$, $x\in X$, and $r>0$ let $B_X(x,r)$ denote the open ball $\{y\in X:d(x,y)<r\}$. Let $\mathbb{Z}_+=\{0,1,\ldots\}$. Our main theorem concerns local well-posedness of the initial-value problem \eqref{Sch1} for small data $s_0\in H^{\sigma_0}_Q(\mathbb{R}^d;\mathbb{S}^2)$, $\sigma_0>(d+1)/2$, $Q\in\mathbb{S}^2$. 

\newtheorem{Main1}{Theorem}[section]
\begin{Main1}\label{Main1}
(a) Assume $\sigma_0>(d+1)/2$ and $Q\in\mathbb{S}^2$. Then there is $\epsilon(\sigma_0)>0$ with the property that for any $s_0\in H^{\infty}_Q(\mathbb{R}^d;\mathbb{S}^2)\cap B_{H^{\sigma_0}_Q(\mathbb{R}^d;\mathbb{S}^2)}(f_Q,\epsilon(\sigma_0))$ there is a unique solution 
\begin{equation*}
s=S^\infty(s_0)\in C([-1,1]:H^{\infty}_Q(\mathbb{R}^d;\mathbb{S}^2))
\end{equation*}
of the initial-value problem \eqref{Sch1}. 

(b) The mapping $s_0\to S^\infty(s_0)$ extends uniquely to a Lipschitz mapping
\begin{equation*}
S^{\sigma_0}:B_{H^{\sigma_0}_Q(\mathbb{R}^d;\mathbb{S}^2)}(f_Q,\epsilon(\sigma_0))\to C([-1,1]:H^{\sigma_0}_Q(\mathbb{R}^d;\mathbb{S}^2)),
\end{equation*}
with the property that $S^{\sigma_0}(s_0)$ is a weak solution of the initial-value problem \eqref{Sch1} for any $s_0\in B_{H^{\sigma_0}_Q(\mathbb{R}^d;\mathbb{S}^2)}(f_Q,\epsilon(\sigma_0))$. 

(c) In addition, for any $\sigma'\in\mathbb{Z}_+$ we have the local Lipschitz bound
\begin{equation}\label{ak1}
\sup_{t\in[-1,1]}d_Q^{\sigma_0+\sigma'}(S^{\sigma_0}(s_0)(t),S^{\sigma_0}(s'_0)(t))\leq C(\sigma_0,\sigma',R)\cdot d_Q^{\sigma_0+\sigma'}(s_0,s'_0)
\end{equation}
for  any $R>0$ and $s_0,s'_0\in B_{H^{\sigma_0}_Q(\mathbb{R}^d;\mathbb{S}^2)}(f_Q,\epsilon(\sigma_0))\cap B_{H^{\sigma_0+\sigma'}_Q(\mathbb{R}^d;\mathbb{S}^2)}(f_Q,R)$. Thus the mapping $S^{\sigma_0}$ restricts to a locally Lipschitz mapping
\begin{equation*}
S^{\sigma_0+\sigma'}:B_{H^{\sigma_0}_Q(\mathbb{R}^d;\mathbb{S}^2)}(f_Q,\epsilon(\sigma_0))\cap H^{\sigma_0+\sigma'}_Q(\mathbb{R}^d;\mathbb{S}^2)\to C([-1,1]:H^{\sigma_0+\sigma'}_Q(\mathbb{R}^d;\mathbb{S}^2)).
\end{equation*}
\end{Main1}

In section \ref{reduc} we use the stereographic projection to reduce Theorem \ref{Main1} to Theorem \ref{Main2}. Then we analyze the resulting derivative Schr\"{o}dinger equation by treating the  nonlinear term as a perturbation. It appears likely that a more careful analysis of the nonlinear interactions, possibly using the ``modified Schr\"{o}dinger map equation'' (cf. \cite{NaStUh} or \cite{ KaKo}), would allow one to extend Theorem \ref{Main1} to the full subcritical range $\sigma_0>d/2$. As in the case of wave maps (for which the regularity theory for small data is much better understood, see \cite{KlMa}, \cite{KlSe}, \cite{Tat1}, \cite{Tat2}, \cite{Ta1}, \cite{Ta2}, \cite{KlRo}, \cite{ShSt}, and \cite{Tat3}), the critical case $\sigma_0=d/2$ is more difficult since, among other things, the critical space $H^{d/2}(\mathbb{R}^d)$ fails to control $L^\infty$. We hope to return to these issues in the future.

The initial-value problem \eqref{Sch1} has been studied extensively (also in the case in which the sphere $\mathbb{S}^2$ is replaced by more general targets). It is known that sufficiently smooth solutions exist locally in time, even for large data (see, for example, \cite{SuSuBa}, \cite{ChShUh}, \cite{DiWa2}, \cite{Ga}, \cite{KePoStTo} and the references therein). Such theorems for (local in time) smooth solutions are proved using variants of the energy method. For low-regularity data, the energy method cannot be applied, and the initial-value problem \eqref{Sch1} has been studied indirectly using the ``modified Schr\"{o}dinger map equation'' (see, for example, \cite{NaStUh}, \cite{NaStUh2}, \cite{KeNa}, and \cite{KaKo}). While existence and uniqueness theorems for this modified Schr\"{o}dinger map equation in certain low-regularity spaces are known (at least in dimension $d=2$), it is not clear whether such theorems can be transfered to the original Schr\"{o}dinger map initial-value problem. Our approach in this paper is more direct, in the sense that we analyze the Schr\"{o}dinger map initial-value problem without passing to the modified Schr\"{o}dinger map equation. As a result of the recursive construction of the solution we obtain a locally Lipschitz flow, which appears to be new even in the case of sufficiently smooth data. Also, our  proof of Theorem \ref{Main1} is self-contained; in particular it does not depend on the existence of smooth solutions.

We describe now some of the ingredients in the proof of Theorem \ref{Main1}. First, using the stereographic projection, we reduce matters to proving Theorem \ref{Main2}. We would then like to analyze the resulting derivative Schr\"{o}dinger equation in some $X^{\sigma ,b}$-type spaces. However, the use of standard $X^{\sigma,b}$ spaces (i.e. spaces defined by suitably weighted norms in the frequency space) seems to lead inevitably to logarithmic divergences, regardless of the amount of smoothness one assumes. To avoid these logarithmic divergences we work with high frequency spaces that have two components: an $X^{\sigma,b}$-type component measured in the frequency space and a normalized $L^{1,2}_{\mathbf{e}}$ (see definition \eqref{vv1}) component measured in the physical space. Such spaces have been used recently in dimension $d=1$ by the authors \cite{IoKe}. The spaces $L^{1,2}_{\mathbf{e}}$ are relevant due to the local smoothing induced by the Schr\"{o}dinger flow. Then we prove suitable linear and nonlinear estimates in these spaces, and conclude Theorem \ref{Main2} using a recursive (perturbative) construction.

The rest of the paper is organized as follows: in section \ref{reduc} we use the stereographic projection to reduce matters to proving Theorem \ref{Main2}. In section \ref{notation} we define our main normed spaces and prove some of their basic properties. In section \ref{estim} we prove several linear and nonlinear estimates. In section \ref{proof} we use these estimates to complete the proof of Theorem \ref{Main2}.   

\section{Preliminary reductions}\label{reduc}

We start now the proof of Theorem \ref{Main1}. By rotation invariance, we may assume 
\begin{equation}\label{al1}
Q=(0,0,1).
\end{equation}
The uniqueness statement in part (a) is straightforward: assume $$s=(s_1,s_2,s_3),s'=(s'_1,s'_2,s'_3)\in C([-1,1]:H^{\infty}_Q(\mathbb{R}^d;\mathbb{S}^2))$$ are solutions of \eqref{Sch1}. Let $q=s'-s$, so
\begin{equation}\label{Sch9}
\begin{cases}
&\partial_tq=(s+q)\times\Delta_x (s+q)-s\times\Delta_x s\,\text{ on }\,\mathbb{R}^d\times[-1,1];\\
&q(0)=0.
\end{cases}
\end{equation}
We multiply \eqref{Sch9} by $q(t)$ and integrate by parts over $\mathbb{R}^d$ to obtain 
\begin{equation}\label{Sch91}
\begin{split}
\frac{1}{2}\partial_t[\|q(t)\|_{L^2}^2]&=\int_{\mathbb{R}^d}[s(t)\times\Delta_x q(t)]\cdot q(t)\,dx\\
&\leq C_s(||q(t)||_{L^2}^2+\sum_{l=1}^d||\partial_{x_l}q(t)||_{L^2}^2).
\end{split}
\end{equation}
Then we apply $\partial_{x_l}$ to \eqref{Sch9}, multiply by $\partial_{x_l}q(t)$, add up over $l=1,\ldots,d$, and integrate by parts over  $\mathbb{R}^d$. The result  is
\begin{equation}\label{Sch92}
\begin{split}
\frac{1}{2}\partial_t[\sum_{l=1}^d\|\partial_{x_l}q(t)\|_{L^2}^2]&=-\int_{\mathbb{R}^d}[q(t)\times\Delta_x s(t)]\cdot \Delta_xq(t)\,dx\\
&\leq C_s(||q(t)||_{L^2}^2+\sum_{l=1}^d||\partial_{x_l}q(t)||_{L^2}^2).
\end{split}
\end{equation}
Using  \eqref{Sch91} and  \eqref{Sch92}, $q\equiv 0$ as desired.

We start now the construction of the solution $s$. Fix $\sigma_0>(d+1)/2$ as in Theorem \ref{Main1}, and $\epsilon_0$ sufficiently small.\footnote{In this section we could have any $\sigma_0>d/2$; for $\sigma_0>(d+1)/2$ the value of $\epsilon_0$ depends only on the dimension $d$.} In view of the Sobolev imbedding theorem, if $f\in B_{H^{\sigma_0}_Q(\mathbb{R}^d;\mathbb{S}^2)}(f_Q,\epsilon_0)$ then $f$ is continuous and takes values in a small neighborhood of $Q$. Thus, for $f=(f_1,f_2,f_3)\in B_{H^{\sigma_0}_Q(\mathbb{R}^d;\mathbb{S}^2)}(f_Q,\epsilon_0)$ we can define
\begin{equation*}
g=L(f)=\frac{f_1+if_2}{1+f_3}.
\end{equation*}
Clearly, $L(f):\mathbb{R}^d\to\mathbb{C}$ is continuous and takes values in a small neighborhood of $0$. For $g\in B_{H^{\sigma_0}}(0,\epsilon_0)$ we define
\begin{equation*}
f=(f_1,f_2,f_3)=\widetilde{L}(g)=\Big(\frac{g+\overline{g}}{1+g\overline{g}},\frac{(-i)(g-\overline{g})}{1+g\overline{g}},\frac{1-g\overline{g}}{1+g\overline{g}}\Big).
\end{equation*}
Clearly, $\widetilde{L}(g):\mathbb{R}^d\to\mathbb{S}^2$ is continuous and takes value in a small neighborhood of $Q$. We have the following estimates:

\newtheorem{Lemmab1}{Lemma}[section]
\begin{Lemmab1}\label{Lemmab1}
(a) If $f\in H^{\infty}_Q(\mathbb{R}^d;\mathbb{S}^2)\cap B_{H^{\sigma_0}_Q(\mathbb{R}^d;\mathbb{S}^2)}(f_Q,\epsilon_0)$ then $L(f)\in H^\infty\cap B_{H^{\sigma_0}}(0,C\epsilon_0)$ and
\begin{equation}\label{al2}
\|L(f)-L(f')\|_{H^\sigma}\leq C(\sigma,d_Q^\sigma(f_Q,f),d_Q^\sigma(f_Q,f')) \cdot d_Q^\sigma(f,f'),
\end{equation}
for any $f,f'\in H^{\infty}_Q(\mathbb{R}^d;\mathbb{S}^2)\cap B_{H^{\sigma_0}_Q(\mathbb{R}^d;\mathbb{S}^2)}(f_Q,\epsilon_0)$ and $\sigma\geq\sigma_0$.

(b) If $g\in H^\infty\cap B_{H^{\sigma_0}}(0,\epsilon_0)$ then $\widetilde{L}(g)\in H^{\infty}_Q(\mathbb{R}^d;\mathbb{S}^2)\cap B_{H^{\sigma_0}_Q(\mathbb{R}^d;\mathbb{S}^2)}(f_Q,C\epsilon_0)$ and 
\begin{equation}\label{al3}
d_Q^\sigma(\widetilde{L}(g),\widetilde{L}(g'))\leq C(\sigma,\|g\|_{H^\sigma},\|g'\|_{H^\sigma}) \cdot \|g-g'\|_{H^\sigma},
\end{equation}
for any $g,g'\in H^\infty\cap B_{H^{\sigma_0}}(0,\epsilon_0)$ and $\sigma\geq\sigma_0$.
\end{Lemmab1}

\begin{proof}[Proof of Lemma \ref{Lemmab1}] In view of the definitions, for both part (a) and part (b) it suffices to prove that if $\sigma\geq\sigma_0$  then
\begin{equation}\label{al4}
\|h\cdot h'\|_{H^\sigma}\leq M(\sigma)\cdot (\|h\|_{H^{\sigma_0}}\cdot \|h'\|_{H^\sigma}+\|h\|_{H^{\sigma}}\cdot \|h'\|_{H^{\sigma_0}})\text{ for any }h,h'\in H^\infty,
\end{equation}
for some increasing function $M:[\sigma_0,\infty)\to[1,\infty)$, and
\begin{equation}\label{al5}
\|h\cdot (1+h')^{-1}\|_{H^\sigma}\leq C(\sigma,\|h'\|_{H^\sigma})\cdot \|h\|_{H^\sigma}\text{ for any }h\in H^\infty,h'\in H^\infty\cap B_{H^{\sigma_0}}(0,\varepsilon_0), 
\end{equation}
provided that $\varepsilon_0$ is sufficiently small. The inequality \eqref{al4} is well-known, using the fact $H^\sigma(\mathbb{R}^d)$ is a Banach algebra for any $\sigma>d/2$ and the Littlewood--Paley product trichotomy (with $M(\sigma)$ depending only on $\sigma$ and $d$). For \eqref{al5} it suffices to prove that 
\begin{equation}\label{al6}
\|h\cdot (h')^n\|_{H^\sigma}\leq 2^{-n}C(\sigma,\|h'\|_{H^\sigma})\cdot \|h\|_{H^\sigma}
\end{equation}
for $h,h'$ as in \eqref{al5}, $\sigma\geq \sigma_0$, and $n\in\mathbb{Z}_+$. The inequality \eqref{al6} clearly holds for $n=0$.

We turn now to the proof  of \eqref{al6} for $n\geq 1$. Let $b_{n,\sigma}=\|h\cdot (h')^n\|_{H^\sigma}$. Using \eqref{al4} we have
\begin{equation}\label{al7}
b_{n,\sigma_0}\leq (M(\sigma_0)||h'||_{H^{\sigma_0}})^n\cdot \|h\|_{H^{\sigma_0}} \text{ for }n=0,1,\ldots,
\end{equation}
which gives \eqref{al6} in the case $\sigma=\sigma_0$. Assume now that $\sigma\leq \sigma_0+2$. Then, using \eqref{al4} and \eqref{al7},
\begin{equation}\label{al8}
\begin{split}
b_{n,\sigma}&\leq M(\sigma_0+2)\cdot ||h'||_{H^{\sigma_0}}\cdot b_{n-1,\sigma}+M(\sigma_0+2)\cdot ||h'||_{H^\sigma}\cdot b_{n-1,\sigma_0}\\
&\leq (1/3)b_{n-1,\sigma}+C(\sigma,||h'||_{H^\sigma})\cdot 2^{-n}||h||_{H^\sigma}.
\end{split}
\end{equation}
Since $b_{0,\sigma}=||h||_{H^\sigma}$, the bound \eqref{al6} follows easily from \eqref{al8} in this case. 

Finally, assume that $\sigma\geq \sigma_0+2$. We may assume that  the bound \eqref{al6} for $\sigma-1$ holds, and use the Kato-Ponce commutator estimate \cite[Lemma XI]{KaPo}: if $f,g\in H^\infty(\mathbb{R}^d)$ and $\sigma>0$ then
\begin{equation*}
||J^\sigma(fg)-fJ^\sigma(g)||_{L^2}\leq C_\sigma(||\nabla f||_{L^\infty}||J^{\sigma-1}g||_{L^2}+||J^\sigma f||_{L^2}||g||_{L^\infty}),
\end{equation*}
where $J^\sigma$ is the operator defined by the multiplier $\xi\to (1+|\xi|^2)^{\sigma/2}$. We apply this  inequality with $f=h'$, $g=h\cdot (h')^{n-1}$. The result is
\begin{equation*}
\begin{split}
b_{n,\sigma}&\leq C||h'||_{H^{\sigma_0}}\cdot b_{n-1,\sigma}+C_\sigma ||h'||_{H^{\sigma_0+1}}\cdot b_{n-1,\sigma-1}+C_\sigma ||h'||_{H^\sigma}||h||_{H^{\sigma_0}}\cdot 2^{-n}\\
&\leq (1/3)b_{n-1,\sigma}+C(\sigma,||h'||_{H^\sigma})\cdot 2^{-n}||h||_{H^\sigma},
\end{split}
\end{equation*}
using the induction hypothesis on $b_{n-1,\sigma-1}$. The bound \eqref{al6} follows in the case $\sigma\geq \sigma_0+2$.
\end{proof}

A direct computation shows that if $u\in C([-1,1]:H^\infty)$ is a solution of the equation
\begin{equation*}
(i\partial_t+\Delta_x)u=\frac{2\overline{u}}{1+u\overline{u}}\sum_{j=1}^d(\partial_{x_j}u)^2\text{ on }\mathbb{R}^d\times[-1,1],
\end{equation*}
then the function $s\in C([-1,1]:H^{\infty}_Q(\mathbb{R}^d;\mathbb{S}^2))$, $s(t)=\widetilde{L}(u(t))$ is a solution of the Schr\"{o}dinger map equation
\begin{equation*}
\partial_ts=s\times\Delta_x s\,\text{ on }\,\mathbb{R}^d\times[-1,1].
\end{equation*}
In view of Lemma \ref{Lemmab1}, for Theorem \ref{Main1} it suffices to prove the following:

\newtheorem{Main2}[Lemmab1]{Theorem}
\begin{Main2}\label{Main2}
(a) Assume $\sigma_0>(d+1)/2$. Then there is $\epsilon(\sigma_0)>0$ with the property that for any $\phi \in H^\infty \cap B_{H^{\sigma_0}}(0,\epsilon(\sigma_0))$ there is a solution 
\begin{equation*}
u=\widetilde{S}^\infty(\phi)\in C([-1,1]:H^\infty)
\end{equation*}
of the initial-value problem
\begin{equation}\label{Sch5}
\begin{cases}
&(i\partial_t+\Delta_x)u=2\overline{u}(1+u\overline{u})^{-1}\sum_{j=1}^d(\partial_{x_j}u)^2\text{ on }\mathbb{R}^d\times[-1,1];\\
&u(0)=\phi.
\end{cases}
\end{equation}

(b) The mapping $\phi\to \widetilde{S}^\infty(\phi)$ extends uniquely to a Lipschitz mapping 
\begin{equation*}
\widetilde{S}^{\sigma_0}:B_{H^{\sigma_0}}(0,\epsilon(\sigma_0))\to C([-1,1]:H^{\sigma_0}),
\end{equation*}
with the property that $\widetilde{S}^{\sigma_0}(\phi)$ is a weak solution of the initial-value problem \eqref{Sch5} for any $\phi\in B_{H^{\sigma_0}}(0,\epsilon(\sigma_0))$.

(c) In addition, for any $\sigma'\in\mathbb{Z}_+$ we have the local Lipchitz bound
\begin{equation}\label{by1}
\sup_{t\in[-1,1]}\|\widetilde{S}^{\sigma_0}(\phi)(t)-\widetilde{S}^{\sigma_0}(\phi')(t)\|_{H^{\sigma_0+\sigma'}}\leq C(\sigma_0,\sigma',R)\cdot \|\phi-\phi'\|_{H^{\sigma_0+\sigma'}},
\end{equation}
for any $R>0$ and $\phi,\phi'\in B_{H^{\sigma_0}}(0,\epsilon(\sigma_0))\cap B_{H^{\sigma_0+\sigma'}}(0,R)$. Thus the mapping $\widetilde{S}^{\sigma_0}$ restricts to a locally Lipschitz mapping
\begin{equation*}
\widetilde{S}^{\sigma_0+\sigma'}:B_{H^{\sigma_0}}(0,\epsilon(\sigma_0))\cap H^{\sigma_0+\sigma'}\to C([-1,1]:H^{\sigma_0+\sigma'}).
\end{equation*} 
\end{Main2}

\section{Notation and preliminary lemmas}\label{notation}

In this section we summarize most of the notation, define our main normed spaces,\footnote{It is likely that only minor changes would be needed to guarantee that all of our normed spaces are in fact Banach spaces. We do not need this, however, since the limiting argument in section \ref{proof} takes place the Banach spaces $C([-1,1]:H^{\sigma})$.}  and prove some of their basic properties. For $l=1,\ldots,d+1$ let $\mathcal{F}_{(l)}$ and $\mathcal{F}_{(l)}^{-1}$ denote the Fourier transform operator and the inverse Fourier transform operator acting on $\mathcal{S}'(\mathbb{R}^l)$. 

For $l=1,\ldots,d$ we fix $\eta_0^{(l)}:\mathbb{R}^l\to[0,1]$ smooth radial functions supported in the sets $\{\xi\in\mathbb{R}^l:|\xi|\leq 8/5]\}$, equal to $1$ in the sets $\{\xi\in\mathbb{R}^l:|\xi|\leq 5/4]\}$, and with the property that
\begin{equation}\label{gu1}
\sum_{k=0}^\infty\eta_k^{(l)}\equiv 1\text{ where }\eta_k^{(l)}(\xi)=\eta_0^{(l)}(\xi/2^k)-\eta_0^{(l)}(\xi/2^{k-1}).
\end{equation}

We define now the normed spaces $X_k$ and $Y_k$. For $k\in\mathbb{Z}_+$ and $j\in\mathbb{Z}_+$ let
\begin{equation*}
\begin{cases}
&D_{k,j}=\{(\xi,\tau)\in\mathbb{R}^d\times\mathbb{R}:|\xi|\in[2^{k-1},2^{k+1}]\text{ and }|\tau+|\xi|^2|\leq 2^{j+1}\}\text{ if }k\geq 1;\\
&D_{k,j}=\{(\xi,\tau)\in\mathbb{R}^d\times\mathbb{R}:|\xi|\leq 2\text{ and }|\tau+|\xi|^2|\leq 2^{j+1}\}\text{ if }k=0.
\end{cases}
\end{equation*}
Let $D_{k,\infty}=\cup_{j\geq 0}D_{k,j}$. We define first the normed spaces
\begin{equation}\label{v1}
\begin{split}
X_k=\{f\in L^2(\mathbb{R}^d\times&\mathbb{R}):f\text{ supported in }D_{k,\infty}\text { and } \\
&\|f\|_{X_k}=\sum_{j=0}^\infty 2^{j/2}\|\eta_j^{(1)}(\tau+|\xi|^2)\cdot f\|_{L^2}<\infty\}.
\end{split}
\end{equation}

The spaces $X_k$ are not sufficient for our estimates, due to various logarithmic divergences. For any vector $\mathbf{e}\in\mathbb{S}^{d-1}$ let $$P_{\mathbf{e}}=\{\xi\in\mathbb{R}^d:\xi\cdot\mathbf{e}=0\}$$ with the induced Euclidean measure. Also, let
\begin{equation*}
D_{k,j}^{\mathbf{e}}=\{(\xi,\tau)\in D_{k,j}:\xi\cdot\mathbf{e}\geq |\xi|/2\}\text{ for }j\in \mathbb{Z}_+\text{ and }j=\infty.
\end{equation*}
For $p,q\in[1,\infty]$ we define the normed spaces $L^{p,q}_{\mathbf{e}}=L^{p,q}_{\mathbf{e}}(\mathbb{R}^d\times\mathbb{R})$,
\begin{equation}\label{vv1}
\begin{split}
L^{p,q}_{\mathbf{e}}&=\{f\in L^2(\mathbb{R}^d\times \mathbb{R}):\\
&\|f\|_{L^{p,q}_{\mathbf{e}}}=\Big[\int_{\mathbb{R}}\Big[\int_{P_\mathbf{e}\times \mathbb{R}}|f(r\mathbf{e}+v,t)|^q\,dvdt\Big]^{p/q}\,dr\Big]^{1/p}<\infty\}.
\end{split}
\end{equation}
Then, for $k\geq 100$ and $\mathbf{e}\in\mathbb{S}^{d-1}$, we define the normed spaces
\begin{equation}\label{v2}
\begin{split}
Y_k^\mathbf{e}=\{f\in L^2&(\mathbb{R}^d\times\mathbb{R}):f\text{ supported in }D_{k,\infty}^\mathbf{e}\text { and } \\
&\|f\|_{Y_k^\mathbf{e}}=2^{-k/2}\|\mathcal{F}_{(d+1)}^{-1}[(\tau+|\xi|^2+i)\cdot f]\|_{L^{1,2}_{\mathbf{e}}}<\infty\}.
\end{split}
\end{equation}
For simplicity of notation, we define $Y_k^\mathbf{e}=\{0\}$ for $k=0,1,\ldots,99$.

We fix $L$ large and $\mathbf{e}_1,\ldots,\mathbf{e}_L\in\mathbb{S}^{d-1}$, $\mathbf{e}_l\neq \mathbf{e}_{l'}$ if $l\neq l'$, with the property that
\begin{equation*}
\text{ for any }\mathbf{e}\in\mathbb{S}^{d-1}\text{ there is }l\in\{1,\ldots,L\}\text{ such that }|\mathbf{e}-\mathbf{e}_l|\leq 2^{-50}.
\end{equation*}
We assume in addition that if $\mathbf{e}\in\{\mathbf{e}_1,\ldots,\mathbf{e}_L\}$ then $-\mathbf{e}\in\{\mathbf{e}_1,\ldots,\mathbf{e}_L\}$. For $k\in\mathbb{Z}_+$ we define
\begin{equation}\label{v3'}
Z_k=X_k+Y_{k}^{\mathbf{e}_1}+\ldots+Y_k^{\mathbf{e}_L}.
\end{equation}

We prove now several estimates. In view of the definitions, if $f\in Z_k$ then we can write
\begin{equation}\label{mi1}
\begin{cases}
&f=\sum\limits_{j\in\mathbb{Z}_+}g_j+f_{\mathbf{e}_1}+\ldots+f_{\mathbf{e}_L}\text{ where }\,g_j\text{ is supported in }D_{k,j},\,f_{\mathbf{e}_l}\in Y_k^{\mathbf{e}_l};\\
&\sum\limits_{j\in\mathbb{Z}_+}2^{j/2}||g_j||_{L^2}+||f_{\mathbf{e}_1}||_{Y_k^{\mathbf{e}_1}}+\ldots+||f_{\mathbf{e}_L}||_{Y_k^{\mathbf{e}_L}}\leq 2\|f\|_{Z_k}.
\end{cases}
\end{equation}
Also, if $m\in L^\infty(\mathbb{R}^d)$, $\mathcal{F}_{(d)}^{-1}(m)\in L^1(\mathbb{R}^d)$, and $f\in Z_k$,  then $m(\xi)\cdot f\in Z_k$ and
\begin{equation}\label{mi2}
||m(\xi)\cdot f||_{Z_k}\leq C||\mathcal{F}_{(d)}^{-1}(m)||_{L^1(\mathbb{R}^d)}\cdot ||f||_{Z_k}.
\end{equation}
For simplicity of notation, for $k\in\mathbb{Z}_+$ and $l\in[0,60]\cap\mathbb{Z}$ we define the smooth functions $\chi_{k,l}:\mathbb{R}\to[0,1]$,
\begin{equation}\label{si100}
\begin{cases}
&\chi_{k,l}(r)=[1-\eta_0^{(1)}(r/2^{k-l})]\cdot \mathbf{1}_{[0,\infty)}(r)\text{ if }k\geq 100;\\
&\chi_{k,l}(r)\equiv 1\text{ if }k\leq 99.
\end{cases}
\end{equation}
We show first that the spaces $Z_k$ are logarithmic modifications of the spaces $X_k$.

\newtheorem{Lemmas1}{Lemma}[section]
\begin{Lemmas1}\label{Lemmas1}
If $k,j\in\mathbb{Z}_+$ and $f\in Z_k$ then
\begin{equation}\label{vvs1}
\|f\cdot \eta_j^{(1)}(\tau+|\xi|^2)\|_{X_k}\leq C\|f\|_{Z_k}\text{ and }\|f\|_{X_k}\leq C(k+1)\|f\|_{Z_k}.
\end{equation} 
\end{Lemmas1}

\begin{proof}[Proof of Lemma \ref{Lemmas1}] Clearly, we may assume $k\geq 100$ and $f=f_{\mathbf{e}}\in Y_k^{\mathbf{e}}$, for some $\mathbf{e}\in\{\mathbf{e}_1,\ldots,\mathbf{e}_L\}$. Let
\begin{equation}\label{gu45}
h_{\mathbf{e}}(x)=2^{-k/2}\mathcal{F}_{(d+1)}^{-1}[(\tau+|\xi|^2+i)\cdot f_{\mathbf{e}}](x).
\end{equation}
Thus
\begin{equation*}
f_{\mathbf{e}}(\xi,\tau)=\mathbf{1}_{D_{k,\infty}^\mathbf{e}}(\xi,\tau)\cdot\frac{2^{k/2}}{\tau+|\xi|^2+i}\mathcal{F}_{(d+1)}(h_{\mathbf{e}})(\xi,\tau).
\end{equation*}
In view of the definitions, for \eqref{vvs1} it suffices to prove that
\begin{equation}\label{gu21}
2^{k/2}2^{-j/2}\|\mathbf{1}_{D_{k,j}}(\xi,\tau)\cdot \chi_{k,40}(\xi\cdot\mathbf{e})\cdot\mathcal{F}_{(d+1)}(h)\|_{L^2_{\xi,\tau}}\leq C(1+2^{j-2k})^{-1/2}\|h\|_{L^{1,2}_{\mathbf{e}}}
\end{equation}
for any $h\in\mathcal{S}(\mathbb{R}^d\times\mathbb{R})$ and $j\in\mathbb{Z}_+$. We write $\xi=\xi_1\mathbf{e}+\xi'$, $x=x_1\mathbf{e}+x'$, $x_1,\xi_1\in\mathbb{R}$, $x',\xi'\in P_{\mathbf{e}}$. Let
\begin{equation*}
h'(x_1,\xi',\tau)=\int_{P_{\mathbf{e}}\times\mathbb{R}}h(x_1\mathbf{e}+x',t)e^{-i(x'\cdot \xi'+t\tau)}\,dx'dt.
\end{equation*}
By Plancherel theorem,
\begin{equation*}
\|h\|_{L^{1,2}_{\mathbf{e}}}=C\|h'\|_{L^1_{x_1}L^2_{\xi',\tau}}.
\end{equation*}
Thus, for \eqref{gu21}, it suffices to prove that
\begin{equation*}
\begin{split}
2^{(k-j)/2}\Big|\Big|\mathbf{1}_{D_{k,j}}(\xi,\tau)\cdot \chi_{k,40}(\xi_1)\int_\mathbb{R}&h'(x_1,\xi',\tau)e^{-ix_1\xi_1}\,dx_1\Big|\Big|_{L^2_{\xi_1,\xi',\tau}}\\\
&\leq C(1+2^{j-2k})^{-1/2}\|h'\|_{L^1_{x_1}L^2_{\xi',\tau}}.
\end{split}
\end{equation*}
This follows easily since for any $(\xi',\tau)\in P_{\mathbf{e}}\times\mathbb{R}$ the measure of the set $\{\xi_1:\xi_1\approx 2^k\text{ and }|\tau+\xi_1^2+|\xi'|^2|\leq 2^{j+1}\}$ is bounded by $C\min(2^{j-k},2^{k})$.
\end{proof}

The proof of Lemma \ref{Lemmas1} shows also that for $k\leq 99$
\begin{equation}\label{gu76}
||f||_{X_k}\leq C2^{-k/2}||\mathcal{F}_{(d+1)}^{-1}[(\tau+|\xi|^2+i)\cdot f]||_{L^{1,2}_{\mathbf{e}}}
\end{equation}
for any $\mathbf{e}\in\mathbb{S}^{d-1}$ and $f$ supported in $D_{k,\infty}$. We prove now a local-smoothing estimate. 

\newtheorem{Lemmas2}[Lemmas1]{Lemma}
\begin{Lemmas2}\label{Lemmas2}
If $k\in\mathbb{Z}_+$, $\mathbf{e}\in\mathbb{S}^{d-1}$, and $f\in Z_k$ then
\begin{equation}\label{gu40}
\|\mathcal{F}_{(d+1)}^{-1}[f\cdot \chi_{k,30}(\xi\cdot\mathbf{e})]\|_{L^{\infty,2}_{\mathbf{e}}}\leq C2^{-k/2}\|f\|_{Z_k}.
\end{equation} 
\end{Lemmas2}
\begin{proof}[Proof of Lemma \ref{Lemmas2}]  We write $\xi=\xi_1\mathbf{e}+\xi'$, $\xi_1\in\mathbb{R}$, $\xi' \in P_{\mathbf{e}}$. Using the Plancherel theorem and the definitions, for \eqref{gu40} it suffices to prove that for any $x_1\in\mathbb{R}$
\begin{equation}\label{gu401}
\Big|\Big|\int_{\mathbb{R}}f(\xi_1\mathbf{e}+\xi',\tau)\chi_{k,30}(\xi_1)e^{ix_1\xi_1}\,d\xi_1\Big|\Big|_{L^2_{\xi',\tau}}\leq C2^{-k/2}\|f\|_{Z_k}.
\end{equation}

We use the representation \eqref{mi1}. Assume first that $f=g_j$. In view of the definitions, it suffices to prove that if $j\geq 0$ and $g_j$ is supported in $D_{k,j}$ then
\begin{equation}\label{gu41}
\Big|\Big|\int_{\mathbb{R}}g_j(\xi_1\mathbf{e}+\xi',\tau)\chi_{k,40}(\xi_1)e^{ix_1\xi_1}\,d\xi_1\Big|\Big|_{L^2_{\xi',\tau}}\leq C2^{-k/2}2^{j/2}\|g_j\|_{L^2}.
\end{equation}
The bound  \eqref{gu41} is a consequence of Plancherel theorem for $k\leq 99$. Assume $k\geq 100$ and let $g_j^\#(\xi_1,\xi',\mu)=g_j(\xi_1\mathbf{e}+\xi',\mu-\xi_1^2-|\xi'|^2)$, so $g_j^\#$ is supported in the set $\{(\xi_1,\xi',\mu)\in\mathbb{R}\times P_\mathbf{e}\times\mathbb{R}:\xi_1^2+|\xi'|^2\in[2^{2k-2},2^{2k+2}],\,|\mu|\leq 2^{j+1}\}$. Using H\"{o}lder's inequality and the fact that $\chi_{k,40}$ is supported in the interval $[c2^{k},\infty)$, the left-hand side of \eqref{gu41} is dominated by
\begin{equation*}
\begin{split}
C&\sup_{||h||_{L^2(P_\mathbf{e}\times\mathbb{R})=1}}\int_{\mathbb{R}\times P_\mathbf{e}\times\mathbb{R}}|g_j(\xi_1\mathbf{e}+\xi',\tau)|\chi_{k,40}(\xi_1)h(\xi',\tau)\,d\xi_1d\xi'd\tau\\
&=C\sup_{||h||_{L^2(P_\mathbf{e}\times\mathbb{R})=1}}\int_{\mathbb{R}\times P_\mathbf{e}\times\mathbb{R}}|g_j^\#(\xi_1,\xi',\mu)|\chi_{k,40}(\xi_1)h(\xi',\mu-\xi_1^2)\,d\xi_1d\xi'd\mu\\
&\leq C2^{-k/2}\int_{\mathbb{R}}\Big[\int_{\mathbb{R}\times P_\mathbf{e}}|g_j^\#(\xi_1,\xi',\mu)|^2\,d\xi_1d\xi' \Big]^{1/2}\,d\mu,
\end{split}
\end{equation*}
which gives \eqref{gu41}

Assume now that $f=f_{\mathbf{e}'}\in Y^{\mathbf{e}'}_k$, $\mathbf{e'}\in\{\mathbf{e}_1,\ldots,\mathbf{e}_L\}$, $k\geq 100$, and define $h_{\mathbf{e}'}$ as in \eqref{gu45}. Notice also that
\begin{equation*}
||f_{\mathbf{e}'}\cdot [1-\eta_0^{(1)}((\tau+|\xi|^2)/2^{2k-100})]||_{X_k}\leq C||f_{\mathbf{e}'}||_{Y_k^{\mathbf{e}'}},
\end{equation*}
using  \eqref{gu21}. Since the inequality \eqref{gu401} was already proved for $f\in X_k$, it suffices to show that for any $x_1\in\mathbb{R}$
\begin{equation}\label{gu47}
\begin{split}
\Big|\Big|\int_{\mathbb{R}}\eta_0^{(1)}((\tau+|\xi|^2)/&2^{2k-100})\,f_{\mathbf{e}'}(\xi_1\mathbf{e}+\xi',\tau)\\
&\times\chi_{k,30}(\xi_1)e^{ix_1\xi_1}\,d\xi_1\Big|\Big|_{L^{2}_{\xi',\tau}}\leq C2^{-k/2}\|f_{\mathbf{e}'}\|_{Y_k^{\mathbf{e}'}}.
\end{split}
\end{equation}

We examine first the support in $(\xi',\tau)$ of the function obtained after taking the $\xi_1$ integral in the left-hand side of \eqref{gu47}. We fix a vector $\mathbf{e}^\perp\in\mathbb{S}^{d-1}\cap P_{\mathbf{e}}$ and a real number $\theta\in[0,2\pi)$ with the property that
\begin{equation}\label{def1}
\mathbf{e}'=\cos\theta\cdot \mathbf{e}+\sin\theta\cdot \mathbf{e}^\perp.
\end{equation}
The choice of $\mathbf{e}^\perp$ and $\theta$ is unique (up to signs) unless $\mathbf{e}'=\pm\mathbf{e}$. The function obtained after taking the $\xi_1$ integral in the left-hand side of \eqref{gu47} is supported in the set 
\begin{equation}\label{gu48}
\begin{split}
S=\{(\xi',\tau)\in P_{\mathbf{e}}\times\mathbb{R}&:-\tau-|\xi'|^2\in[2^{2k-80},2^{2k+10}], |\xi'|\leq 2^{k+1},\text{ and }\\
&(-\tau-|\xi'|^2)^{1/2}\cos\theta+(\xi'\cdot\mathbf{e}^{\perp})\sin\theta\geq 2^{k-10}\},
\end{split}
\end{equation}
and  the integral in $\xi_1$ is taken over the interval $\xi_1\in[2^{k-30},2^{k+1}]$. As in the proof of Lemma \ref{Lemmas1}, we can write
\begin{equation*}
f_{\mathbf{e}'}(\xi,\tau)=\chi_{k,5}(\xi\cdot \mathbf{e}')\cdot\frac{2^{k/2}}{\tau+|\xi|^2+i}\mathcal{F}_{(d+1)}(h_{\mathbf{e}'})(\xi,\tau)\text{ where }\|f_{\mathbf{e}'}\|_{Y_k^{\mathbf{e}'}}=C\|h_{\mathbf{e}'}\|_{L^{1,2}_{\mathbf{e}'}}.
\end{equation*}
Thus, for \eqref{gu47} it suffices to prove that
\begin{equation}\label{gu47'}
\begin{split}
\Big|\Big|&\mathbf{1}_S(\xi',\tau)\int_{\mathbb{R}}\eta_0^{(1)}((\tau+|\xi|^2)/2^{2k-100})\mathcal{F}_{(d+1)}(h_{\mathbf{e}'})(\xi_1\mathbf{e}+\xi',\tau)\\
&\times\chi_{k,5}(\xi\cdot \mathbf{e}')(\tau+|\xi|^2+i)^{-1}\chi_{k,30}(\xi_1)e^{ix_1\xi_1}\,d\xi_1\Big|\Big|_{L^{2}_{\xi',\tau}}\leq C2^{-k}\|h_{\mathbf{e}'}\|_{L^{1,2}_{\mathbf{e}'}}.
\end{split}
\end{equation}

Let $M=M(\xi',\tau)=(-\tau-|\xi'|^2)^{1/2}$. Elementary estimates using the definitions show that we can approximate
\begin{equation}\label{gu49}
\begin{split}
&\mathbf{1}_{S}(\xi',\tau)\cdot\eta_0^{(1)}((\xi_1^2-M^2)/2^{2k-100})\cdot (\xi_1^2-M^2+i)^{-1}\cdot \chi_{k,5}(\xi\cdot \mathbf{e}')\cdot \chi_{k,30}(\xi_1)\\
&=\mathbf{1}_{S}(\xi',\tau)\cdot \chi_{k,30}(M)\cdot \chi_{k,5}((M\mathbf{e}+\xi')\cdot\mathbf{e}')\frac{\eta_0^{(1)}((\xi_1-M)/2^{k-100})}{2M\cdot(\xi_1-M+i/2^{k})}+E
\end{split}
\end{equation}
where
\begin{equation}\label{gu50}
\begin{split}
|E(\xi,\tau)|&\leq C\cdot\mathbf{1}_{S}(\xi',\tau)\cdot\mathbf{1}_{[2^{k-35},2^{k+10}]}(\xi_1)\cdot\mathbf{1}_{[2^{k-35},2^{k+10}]}(\xi\cdot\mathbf{e}')\\
&\times\eta_0^{(1)}((\xi_1^2-M^2)/2^{2k+10})\cdot\big[2^{-2k}+\big(1+|\xi_1^2-M^2|\big)^{-2}\big].
\end{split}
\end{equation}

We substitute the identity \eqref{gu49} into \eqref{gu47'}. We handle first the error term: using \eqref{gu21} and \eqref{gu41}
\begin{equation*}
\begin{split}
&\Big|\Big|\int_{\mathbb{R}}E(\xi,\tau)\mathcal{F}_{(d+1)}(h_{\mathbf{e}'})(\xi_1\mathbf{e}+\xi',\tau)e^{ix_1\xi_1}\,d\xi_1\Big|\Big|_{L^{2}_{\xi',\tau}}\\
&\leq \sum_{j\leq 2k+C}\Big|\Big|\int_{\mathbb{R}}\eta_j^{(1)}(\tau+|\xi|^2)\cdot E(\xi,\tau)\mathcal{F}_{(d+1)}(h_{\mathbf{e}'})(\xi_1\mathbf{e}+\xi',\tau)e^{ix_1\xi_1}\,d\xi_1\Big|\Big|_{L^{2}_{\xi',\tau}}\\
&\leq C\sum_{j\leq 2k+C}2^{-(k-j)/2}||\eta_j^{(1)}(\tau+|\xi|^2)\cdot E(\xi,\tau)\cdot \mathcal{F}_{(d+1)}(h_{\mathbf{e}'})(\xi,\tau)||_{L^{2}_{\xi,\tau}}\\
&\leq C\sum_{j\leq 2k+C}2^{-(k-j)/2}\cdot(2^{-2k}+2^{-2j})\cdot 2^{-(k-j)/2}\|h_{\mathbf{e}'}\|_{L^{1,2}_{\mathbf{e}'}}\\
&\leq C2^{-k}\|h_{\mathbf{e}'}\|_{L^{1,2}_{\mathbf{e}'}},
\end{split}
\end{equation*}
which agrees with \eqref{gu47'}.

We estimate now the contribution of the first term in the right-hand side of \eqref{gu49}. Since $M(\xi_2,\tau)\approx 2^{k}$ in $S$, it suffices to prove that
\begin{equation}\label{gu51}
\begin{split}
\Big|\Big|\mathbf{1}_S(\xi',\tau)&\int_{\mathbb{R}}\frac{\eta_0^{(1)}((\xi_1-M(\xi',\tau))/2^{k-100})}{\xi_1-M(\xi',\tau)+i/2^{k}}\\
&\times\mathcal{F}_{(d+1)}(h)(\xi,\tau)\cdot e^{ix_1\xi_1}\,d\xi_1\Big|\Big|_{L^{2}_{\xi',\tau}}\leq C\|h\|_{L^{1,2}_{\mathbf{e}'}},
\end{split}
\end{equation}
for any $h\in\mathcal{S}(\mathbb{R}^d\times\mathbb{R})$ and $x_1\in\mathbb{R}$. With $\mathbf{e}^\perp$ as $\theta$ as in \eqref{def1}, let
\begin{equation}\label{def2}
{\mathbf{e}'}^\perp=-\sin\theta\cdot \mathbf{e}+\cos\theta\cdot \mathbf{e}^\perp.
\end{equation}
Let
\begin{equation*}
\widetilde{P}_{\mathbf{e},\mathbf{e}^\perp}=\{\xi\in\mathbb{R}^d:\xi\cdot\mathbf{e}=\xi\cdot\mathbf{e}^\perp=0\}=\widetilde{P}_{\mathbf{e}',{\mathbf{e}'}^\perp},
\end{equation*}
and write $\xi=\xi_1\mathbf{e}+\xi_2\mathbf{e}^\perp+\xi''$, $y=y_1\mathbf{e}'+y_2{\mathbf{e}'}^\perp+y''$, $y_1,y_2,\xi_1,\xi_2\in\mathbb{R}$, $y'',\xi''\in\widetilde{P}_{\mathbf{e},\mathbf{e}^\perp}$. For $\tau,r\in\mathbb{R}$ and $\xi''\in\widetilde{P}_{\mathbf{e},\mathbf{e}^\perp}$ let
\begin{equation*}
h'(y_1,r,\xi'',\tau)=\int_{\mathbb{R}\times\widetilde{P}_{\mathbf{e},\mathbf{e}^\perp}\times\mathbb{R}}h(y_1\mathbf{e'}+y_2{\mathbf{e}'}^\perp+y'',t)e^{-i y_2r}e^{-iy''\cdot\xi''}e^{-it\tau}\, dy_2 dy''dt.
\end{equation*}
By Plancherel theorem,
\begin{equation*}
\|h\|_{L^{1,2}_{\mathbf{e}'}}=C||h'||_{L^1_{y_1}L^2_{r,\xi'',\tau}}.
\end{equation*}
Also, using  \eqref{def1} and \eqref{def2},
\begin{equation*}
\begin{split}
\mathcal{F}_{(d+1)}(h)&(\xi_1\mathbf{e}+\xi_2\mathbf{e}^\perp+\xi'',\tau)\\
&=\int_{\mathbb{R}}h'(y_1,-\sin\theta\cdot\xi_1+\cos\theta\cdot\xi_2,\xi'',\tau)e^{-iy_1(\cos\theta\cdot\xi_1+\sin\theta\cdot\xi_2)}\,dy_1.
\end{split}
\end{equation*}
Thus, for \eqref{gu51} it suffices to prove that
\begin{equation}\label{gu52}
\begin{split}
\Big|\Big|\mathbf{1}_S&(\xi_2\mathbf{e}^\perp+\xi'',\tau)\int_{\mathbb{R}}\frac{\eta_0^{(1)}[(\xi_1-M(\xi_2\mathbf{e}^\perp+\xi'',\tau))/2^{k-100}]}{\xi_1-M(\xi_2\mathbf{e}^\perp+\xi'',\tau)+i/2^{k}}\\
&\times h''(-\sin\theta\cdot\xi_1+\cos\theta\cdot\xi_2,\xi'',\tau)\cdot e^{ix_1\xi_1}\,d\xi_1\Big|\Big|_{L^{2}_{\xi_2,\xi'',\tau}}\leq C\|h''\|_{L^2},
\end{split}
\end{equation}
for any compactly supported function $h'':\mathbb{R}\times P_{\mathbf{e},\mathbf{e}^\perp}\times\mathbb{R}\to\mathbb{C}$. Let
\begin{equation*}
h^\ast_{x_1}(r,\xi'',\tau)=\int_{\mathbb{R}}\frac{\eta_0^{(1)}(v/2^{k-100})}{v+i/2^{k}}h''(-\sin\theta\cdot v+r)\cdot e^{ix_1v}\,dv
\end{equation*}
Using the boundedness of the Hilbert transform on $L^2(\mathbb{R})$, $||h^\ast_{x_1}||_{L^2}\leq C||h''||_{L^2}$. Thus, for \eqref{gu52} it suffices to prove that
\begin{equation*}
\begin{split}
||&\mathbf{1}_S(\xi_2\mathbf{e}^\perp+\xi'',\tau)\\
&\times h^\ast_{x_1}(-\sin\theta\cdot M(\xi_2\mathbf{e}^\perp+\xi'',\tau)+\cos\theta\cdot\xi_2,\xi'',\tau)||_{L^{2}_{\xi_2,\xi'',\tau}}\leq C\|h^\ast_{x_1}\|_{L^2}.
\end{split}
\end{equation*}
This follows easily by a change of variables, using the definition \eqref{gu48} of the set $S$. This completes the proof of the lemma.
\end{proof}

We remark  that the proof also gives the following weaker inequality: if $k\in\mathbb{Z}_+$, $\mathbf{e}\in\mathbb{S}^{d-1}$, and $f\in Z_k$ then
\begin{equation}\label{gu90}
\|\mathcal{F}_{(d+1)}^{-1}(f)\|_{L^{\infty,2}_{\mathbf{e}}}\leq C(k+1)\|f\|_{Z_k}.
\end{equation}
For this, using  Lemma \ref{Lemmaa2}, it suffice to prove that
\begin{equation*}
\|\mathcal{F}_{(d+1)}^{-1}(g)\|_{L^{\infty,2}_{\mathbf{e}}}\leq C\|g\|_{X_k}
\end{equation*}
for any $g\in X_k$. We decompose $g=\sum_{j=0}^\infty g_j$, $g_j$ supported in $D_{k,j}$, write $\xi=\xi_1\mathbf{e}+\xi'$, $\xi_1\in\mathbb{R}$, $\xi'\in P_\mathbf{e}$, and  use Plancherel theorem. It remains to prove that
\begin{equation*}
\Big|\Big|\int_{\mathbb{R}}g_j(\xi_1\mathbf{e}+\xi',\tau)e^{ix_1\xi_1}\,d\xi_1\Big|\Big|_{L^2_{\xi',\tau}}\leq C2^{j/2}\|g_j\|_{L^2}\text{ for any  }x_1\in\mathbb{R}.
\end{equation*}
We decompose $g_j=\eta_0^{(1)}(\xi_1)\cdot g_j+(1-\eta_0^{(1)}(\xi_1))\cdot g_j$, and apply H\"older's inequality for the first part and the same argument as in the proof of  \eqref{gu41} for the second part. This completes the proof of \eqref{gu90}.

We will also need a maximal function estimate.

\newtheorem{Lemmaa1}[Lemmas1]{Lemma}
\begin{Lemmaa1}\label{Lemmaa1}
If $k\geq 0$, $f\in Z_k$, and $\mathbf{e}\in\mathbb{S}^{d-1}$ then
\begin{equation}\label{pr40}
||\mathbf{1}_{[-2,2]}(t)\cdot \mathcal{F}_{(d+1)}^{-1}(f)||_{L^{2,\infty}_{\mathbf{e}}}\leq
C2^{(d-1)k/2}(k+1)^2\cdot \|f\|_{Z_k}.
\end{equation}
\end{Lemmaa1}

\begin{proof}[Proof of Lemma \ref{Lemmaa1}] In view of Lemma \ref{Lemmas1}, we may assume $f\in X_k$. Using \eqref{mi1} it suffices to prove that
\begin{equation}\label{pr41}
||\mathbf{1}_{[-2,2]}(t)\cdot \mathcal{F}_{(d+1)}^{-1}(g_j)||_{L^{2,\infty}_{\mathbf{e}}}\leq
C2^{(d-1)k/2}(k+1)\cdot 2^{j/2}\|g_j\|_{L^2}
\end{equation}
for any function $g_j$ supported in $D_{k,j}$. We define $g_j^\#(\xi,\mu)=g_j(\xi,\mu-|\xi|^2)$. The left-hand side of \eqref{pr41} is dominated by
\begin{equation*}
\begin{split}
\int_{[-2^{j+1},2^{j+1}]}\Big|\Big|\mathbf{1}_{[-2,2]}(t)\cdot \int_{\mathbb{R}^d}g^\#_{j}(\xi,\mu)e^{ix\cdot \xi}e^{-it|\xi|^2}\,d\xi\Big|\Big|_{L^{2,\infty}_{\mathbf{e}}}\,d\mu.
\end{split}
\end{equation*}
Thus, for \eqref{pr41} it suffices to prove that
\begin{equation}\label{ar400}
\Big|\Big|\mathbf{1}_{[-2,2]}(t)\int_{\mathbb{R}^d}h(\xi)e^{ix\xi}e^{-it|\xi|^2}\,d\xi\Big|\Big|_{L^{2,\infty}_{\mathbf{e}}}\leq C2^{(d-1)k/2}(k+1)\cdot ||h||_{L^2_\xi},
\end{equation}
for any function $h$ supported in the set $\{\xi\in\mathbb{R}^d:|\xi|\leq 2^{k+1}\}$.

To prove \eqref{ar400}, using a standard $TT^\ast$ argument, it suffices to show that
\begin{equation}\label{ar401}
\begin{split}
\Big|&\Big|\mathbf{1}_{[-4,4]}(t)\int_{\mathbb{R}^{d-1}\times\mathbb{R}}e^{ix_1\xi_1}e^{ix'\cdot \xi'}e^{-it(\xi_1^2+|\xi'|^2)}\\
&\times \eta_0^{(1)}(\xi_1/2^{k+1})\cdot \eta_0^{(d-1)}(\xi'/2^{k+1})\,d\xi_1d\xi'\Big|\Big|_{L^1_{x_1}L^\infty_{x',t}}\leq C2^{(d-1)k}(k+1)^2.
\end{split}
\end{equation}
By stationary phase, for any $\xi'\in\mathbb{R}^{d-1}$
\begin{equation*}
\Big|\int_{\mathbb{R}^{d-1}}e^{ix'\cdot \xi'}e^{-it|\xi'|^2}\eta_0^{(d-1)}(\xi'/2^{k+1})\,d\xi'\Big|\leq C\min(2^{(d-1)k},|t|^{-(d-1)/2}),
\end{equation*}
and
\begin{equation*}
\Big|\int_{\mathbb{R}}e^{ix_1\cdot \xi_1}e^{-it\xi_1^2}\eta_0^{(1)}(\xi_1/2^{k+1})\,d\xi_1\Big|\leq C\min(2^{k},|t|^{-1/2}).
\end{equation*}
In addition, by integration by parts, if $|x_1|\geq 2^{k+10}|t|$ then
\begin{equation*}
\Big|\int_{\mathbb{R}}e^{ix_1\cdot \xi_1}e^{-it\xi_1^2}\eta_0^{(1)}(\xi_1/2^{k+1})\,d\xi_1\Big|\leq C2^k(1+2^k|x_1| )^{-2}.
\end{equation*}
Let $K(x_1,x',t)$ denote the function in the left-hand side of \eqref{ar401}. In view of the three bounds above,
\begin{equation*}
\sup_{|t|\leq 4,\,x'\in\mathbb{R}^{d-1}}|K(x_1,x',t)|\leq C2^{dk}(1+2^k|x_1| )^{-2}+C2^{dk/2}|x_1|^{-d/2}\cdot\mathbf{1}_{[2^{-k},2^k]}( |x_1| ).
\end{equation*}
The bound \eqref{ar401} follows since $d\geq 2$.
\end{proof}

We conclude this section with $L^\infty_tL^2_x$ and $L^\infty_{x,t}$ estimates.

\newtheorem{Lemmaa2}[Lemmas1]{Lemma}
\begin{Lemmaa2}\label{Lemmaa2}
If $k\geq 0$, $t\in\mathbb{R}$, and $f\in Z_k$ then
\begin{equation}\label{lb4}
\sup_{t\in\mathbb{R}}\|\mathcal{F}_{(d+1)}^{-1}(f)(.,t)\|_{L^2_x}
\leq C\|f\|_{Z_k}.
\end{equation}
Thus
\begin{equation}\label{lb44}
\|\mathcal{F}_{(d+1)}^{-1}(f)\|_{L^\infty_{x,t}}
\leq C2^{dk/2}\|f\|_{Z_k}.
\end{equation}
\end{Lemmaa2}

\begin{proof}[Proof of Lemma \ref{Lemmaa2}] By Plancherel theorem it suffices to prove that
\begin{equation}\label{si1}
\Big|\Big|\int_{\mathbb{R}}f(\xi,\tau)e^{it\tau}\,d\tau\Big|\Big|_{L^2_\xi}\leq C\|f\|_{Z_k}
\end{equation}
We use the representation \eqref{mi1}. Assume first that $f=g_j$. Then
\begin{equation*}
\begin{split}
\Big|\Big|\int_\mathbb{R}g_j(\xi,\tau)e^{it\tau}\,d\tau\Big|\Big|_{L^2_\xi}\leq C||g_j(\xi,\tau)||_{L^2_\xi L^1_\tau}\leq C2^{j/2}||g_j||_{L^2_{\xi,\tau}},
\end{split}
\end{equation*}
which proves \eqref{si1} in this case.

Assume now that $k\geq 100$ and $f=f_{\mathbf{e}}\in Y_k^\mathbf{e}$, $\mathbf{e}\in\{\mathbf{e}_1,\ldots,\mathbf{e}_L\}$. We have to prove that
\begin{equation}\label{si1'}
\Big|\Big|\int_{\mathbb{R}}f_{\mathbf{e}}(\xi,\tau)e^{it\tau}\,d\tau\Big|\Big|_{L^2_\xi}\leq C\|f_\mathbf{e}\|_{Y_k^\mathbf{e}}
\end{equation}
We define $h_\mathbf{e}$ as  in \eqref{gu45}, so
\begin{equation*}
f_{\mathbf{e}}(\xi,\tau)=\chi_{k,10}(\xi\cdot\mathbf{e})\cdot\frac{2^{k/2}}{\tau+|\xi|^2+i}\mathcal{F}_{(d+1)}(h_{\mathbf{e}})(\xi,\tau),
\end{equation*}
with $\chi_{k,10}$ as in \eqref{si100}. We write $\xi=\xi_1\mathbf{e}+\xi'$, $x=x_1\mathbf{e}+x'$, $x_1,\xi_1\in\mathbb{R}$, $x',\xi'\in P_{\mathbf{e}}$. For  \eqref{si1'} it suffices  to prove that
\begin{equation}\label{si2}
2^{k/2}\Big|\Big|\chi_{k,10}(\xi_1)\int_{\mathbb{R}}\frac{1}{\tau+|\xi|^2+i}\cdot \mathcal{F}_{(d+1)}(h)(\xi_1\mathbf{e}+\xi',\tau)e^{it\tau}\,d\tau\Big|\Big|_{L^2_\xi}\leq C||h||_{L^{1,2}_{\mathbf{e}}},
\end{equation}
for any $h\in\mathcal{S}(\mathbb{R}^d\times\mathbb{R})$ and $t\in\mathbb{R}$. As in the proof of Lemma  \ref{Lemmas1}, we define
\begin{equation*}
h'(x_1,\xi',\tau)=\int_{P_{\mathbf{e}}\times\mathbb{R}}h(x_1\mathbf{e}+x',t)e^{-i(x'\cdot \xi'+t\tau)}\,dx'dt,
\end{equation*}
so
\begin{equation*}
\mathcal{F}_{(d+1)}(h)(\xi_1\mathbf{e}+\xi',\tau)=\int_{\mathbb{R}}h'(x_1,\xi',\tau)e^{-ix_1\xi_1}\,dx_1\text{ and }\|h\|_{L^{1,2}_{\mathbf{e}}}=C\|h'\|_{L^1_{x_1}L^2_{\xi',\tau}}.
\end{equation*}
Let
\begin{equation*}
h^{\ast}_t(x_1,\xi',\mu)=\int_{\mathbb{R}}\frac{1}{\tau+\mu+i}h'(x_1,\xi',\tau)e^{it\tau}\,d\tau.
\end{equation*}
In view of the boundedness of the Hilbert transform on $L^2(\mathbb{R})$,
\begin{equation*}
||h^{\ast}_t(x_1,\xi',\mu)||_{L^2_{\xi',\mu}}\leq C||h'(x_1,\xi',\tau)||_{L^2_{\xi',\tau}}\text{ for any }x_1,t\in\mathbb{R}.
\end{equation*}
Thus, for \eqref{si2}, it suffices to prove that
\begin{equation*}
2^{k/2}\Big|\Big|\chi_{k,10}(\xi_1)\int_{\mathbb{R}}h^\ast_t(x_1,\xi',|\xi|^2)e^{-ix_1\xi_1}\,dx_1\Big|\Big|_{L^2_\xi}\leq C||h^\ast_t||_{L^1_{x_1}L^2_{\xi',\mu}}.
\end{equation*}
This follows easily by changes of variables. 
\end{proof}
 
\section{Linear and nonlinear estimates}\label{estim}

For $\sigma\geq 0$ we define the normed spaces
\begin{equation}\label{no5}
F^\sigma=\{u\in C(\mathbb{R}:H^\infty):\|u\|_{F^\sigma}^2=\sum_{k=0}^\infty 2^{2\sigma k}\|\eta_k^{(d)}(\xi)\cdot \mathcal{F}_{(d+1)}u\|_{Z_k}^2<\infty\},
\end{equation}
and
\begin{equation}\label{no6}
\begin{split}
N^\sigma=&\{u\in C(\mathbb{R}:H^\infty):\\
&\|u\|_{N^\sigma}^2=\sum_{k=0}^\infty 2^{2\sigma k}\|\eta_k^{(d)}(\xi)\cdot (\tau+|\xi|^2+i)^{-1}\cdot \mathcal{F}_{(d+1)}u\|_{Z_k}^2<\infty\}.
\end{split}
\end{equation}

For $\phi\in H^\infty$ let $W(t)\phi\in
C(\mathbb{R}:H^\infty)$ denote the solution of the free Schr\"{o}dinger evolution
\begin{equation}\label{ni1}
[W(t)\phi](x,t)=c_0\int_{\mathbb{R}^d}e^{ix\cdot\xi}e^{-it|\xi|^2}\mathcal{F}_{(d)}(\phi)(\xi)\,d\xi.
\end{equation}
Assume
$\psi:\mathbb{R}\to[0,1]$ is an even smooth function supported in
the interval $[-8/5,8/5]$ and equal to $1$ in the interval
$[-5/4,5/4]$. We prove first two linear estimates.

\newtheorem{Lemmaq1}{Lemma}[section]
\begin{Lemmaq1}\label{Lemmaq1}
If $\sigma\geq 0$ and $\phi\in H^{\infty}$ then $\psi(t)\cdot [W(t)\phi]\in F^\sigma$ and
\begin{equation*}
\|\psi(t)\cdot [W(t)\phi]\|_{F^{\sigma}}\leq C_\sigma\|\phi\|_{H^\sigma}.
\end{equation*}
\end{Lemmaq1}

\begin{proof}[Proof of Lemma \ref{Lemmaq1}] A straightforward computation shows that
\begin{equation*}
\mathcal{F}_{(d+1)}[\psi(t)\cdot (W(t)\phi)](\xi,\tau)=
\mathcal{F}_{(d)}(\phi)(\xi)\cdot \mathcal{F}_{(1)}(\psi)(\tau+|\xi|^2).
\end{equation*}
Then, directly from the definitions,
\begin{equation*}
\begin{split}
\|\psi(t)\cdot [W(t)\phi]\|_{F^{\sigma}}^2&=\sum_{k\in\mathbb{Z}_+}2^{2\sigma k}\|\eta_k^{(d)}(\xi)\cdot \mathcal{F}_{(d)}(\phi)(\xi)\cdot \mathcal{F}_{(1)}(\psi)(\tau+|\xi|^2)\|_{Z_k}^2\\
&\leq \sum_{k\in\mathbb{Z}_+}2^{2\sigma k}\|\eta_k^{(d)}(\xi)\cdot \mathcal{F}_{(d)}(\phi)(\xi)\cdot \mathcal{F}_{(1)}(\psi)(\tau+|\xi|^2)\|_{X_k}^2\\
&\leq C\sum_{k\in\mathbb{Z}_+}2^{2\sigma k}\|\eta_k^{(d)}(\xi)\cdot \mathcal{F}_{(d)}(\phi)(\xi)\|_{L^2}^2\\
&\leq C_\sigma\|\phi\|_{H^\sigma},
\end{split}
\end{equation*}
as desired.
\end{proof}

\newtheorem{Lemmaq3}[Lemmaq1]{Lemma}
\begin{Lemmaq3}\label{Lemmaq3}
If $\sigma\geq 0$ and $u\in N^{\sigma}$ then $\psi(t)\cdot \int_0^tW(t-s)(u(s))\,ds\in F^\sigma$ and
\begin{equation*}
\Big|\Big|\psi(t)\cdot \int_0^tW(t-s)(u(s))\,ds\Big|\Big|_{F^{\sigma}}\leq C||u||_{N^{\sigma}}.
\end{equation*}
\end{Lemmaq3}

\begin{proof}[Proof of Lemma \ref{Lemmaq3}] A straightforward computation shows that
\begin{equation*}
\begin{split}
\mathcal{F}_{(d+1)}\Big[\psi(t)\cdot& \int_0^tW(t-s)(u(s))ds\Big](\xi,\tau)=\\
&c\int_\mathbb{R}\mathcal{F}_{(d+1)}(u)(\xi,\tau')\frac{\widehat{\psi}(\tau-\tau')-\widehat{\psi}(\tau+|\xi|^2)}{\tau'+|\xi|^2}d\tau',
\end{split}
\end{equation*}
where, for simplicity of notation, $\widehat{\psi}=\mathcal{F}_{(1)}(\psi)$. For $k\in\mathbb{Z}_+$ let
$$f_k(\xi,\tau')=\mathcal{F}_{(d+1)}(u)(\xi,\tau')\cdot \eta_k^{(d)}(\xi)\cdot (\tau'+|\xi|^2+i)^{-1}.$$
For $f\in Z_k$ let
\begin{equation}\label{ar202}
T(f)(\xi,\tau)=\int_\mathbb{R}f(\xi,\tau')\frac{\widehat{\psi}(\tau-\tau')-
\widehat{\psi}(\tau+|\xi|^2)}{\tau'+|\xi|^2}(\tau'+|\xi|^2+i)\,d\tau'.
\end{equation}
In view of the definitions, it suffices to prove that
\begin{equation}\label{ni5}
||T||_{Z_k\to Z_k}\leq C\text{ uniformly in }k\in\mathbb{Z}_+.
\end{equation}

To prove \eqref{ni5} we use the representation \eqref{mi1}. Assume first that $f=g_j$ is supported in $D_{k,j}$. Let $g_j^\#(\xi,\mu')=g_j(\xi,\mu'-|\xi|^2)$ and $[T(g)]^\#(\xi,\mu)=T(g)(\xi,\mu-|\xi|^2)$. Then,
\begin{equation}\label{ni6}
[T(g)]^\#(\xi,\mu)=\int_\mathbb{R}g_j^\#(\xi,\mu')\frac{\widehat{\psi}(\mu-\mu')-
\widehat{\psi}(\mu)}{\mu'}(\mu'+i)\,d\mu'.
\end{equation}
We use the elementary bound
\begin{equation*}
\Big|\frac{\widehat{\psi}(\mu-\mu')-\widehat{\psi}(\mu)}{\mu'}(\mu'+i)\Big|\leq C[(1+|\mu|)^{-4}+(1+|\mu-\mu'|)^{-4}].
\end{equation*}
Then, using \eqref{ni6},
\begin{equation*}
\begin{split}
|T(g)^\#(\xi,\mu)|&\leq C(1+|\mu|)^{-4}\cdot 2^{j/2}\Big[\int_{\mathbb{R}}|g_j^\#(\xi,\mu')|^2\,d\mu'\Big]^{1/2}\\
&+C\mathbf{1}_{[-2^{j+10},2^{j+10}]}(\mu)\int_{\mathbb{R}}|g_j^\#(\xi,\mu')|(1+|\mu-\mu'|)^{-4}\,d\mu'.
\end{split}
\end{equation*}
It follows from the definition of the spaces $X_k$ that
\begin{equation}\label{ni7}
||T||_{X_k\to X_k}\leq C\text{ uniformly in }k\in\mathbb{Z}_+,
\end{equation}
as desired.

Assume now that $f=f_{\mathbf{e}}\in Y^\mathbf{e}_k$,  $k\geq 100$, $\mathbf{e}\in\{\mathbf{e}_1,\ldots,\mathbf{e}_L\}$. We write
\begin{equation*}
f_{\mathbf{e}}(\xi,\tau')=\frac{\tau'+|\xi|^2}{\tau'+|\xi|^2+i}f_{\mathbf{e}}(\xi,\tau')+\frac{i}{\tau'+|\xi|^2+i}f_{\mathbf{e}}(\xi,\tau').
\end{equation*}
Using Lemma \ref{Lemmas1}, $||i(\tau'+|\xi|^2+i)^{-1}f_{\mathbf{e}}(\xi,\tau')||_{X_k}\leq C||g_k||_{Y_k^\mathbf{e}}$. In view of \eqref{ar202} and \eqref{ni7}, for \eqref{ni5} it suffices to prove that
\begin{equation}\label{ni8}
\Big|\Big|\int_\mathbb{R}f_{\mathbf{e}}(\xi,\tau')\widehat{\psi}(\tau-\tau')\,d\tau'\Big|\Big|_{Z_k}+
\Big|\Big|\widehat{\psi}(\tau+|\xi|^2)\int_\mathbb{R}f_{\mathbf{e}}(\xi,\tau')\,d\tau'\Big|\Big|_{X_k}\leq C||f_{\mathbf{e}}||_{Y_k^{\mathbf{e}}}.
\end{equation}
The bound for the second term in the left-hand side of \eqref{ni8} follows from \eqref{si1'} with $t=0$. To bound the first term we write
\begin{equation*}
f_{\mathbf{e}}(\xi,\tau')=f_{\mathbf{e}}(\xi,\tau')\Big[\frac{\tau'+|\xi|^2+i}{\tau+|\xi|^2+i}+\frac{\tau-\tau'}{\tau+|\xi|^2+i}\Big].
\end{equation*}
The first term in the left-hand side of \eqref{ni8} is dominated by
\begin{equation}\label{ni9}
\begin{split}
&C\Big|\Big|(\tau+|\xi|^2+i)^{-1}\int_\mathbb{R}f_{\mathbf{e}}(\xi,\tau')(\tau'+|\xi|^2+i)\widehat{\psi}(\tau-\tau')\,d\tau'\Big|\Big|_{Y_k^{\mathbf{e}}}\\
&+C\Big|\Big|(\tau+|\xi|^2+i)^{-1}\int_\mathbb{R}f_{\mathbf{e}}(\xi,\tau')\widehat{\psi}(\tau-\tau')\cdot (\tau-\tau')\,d\tau'\Big|\Big|_{X_k}.
\end{split}
\end{equation}
For the first term in \eqref{ni9} we use the definition to bound it by $C||f_{\mathbf{e}}||_{Y_k^{\mathbf{e}}}$. For the second term in \eqref{ni9}, it follows from Lemma \ref{Lemmas1} that $\|f_{\mathbf{e}}\|_{L^2_{\xi,\tau'}}\leq C||f_{\mathbf{e}}||_{Y_k^{\mathbf{e}}}$, thus $$
\Big|\Big|\int_\mathbb{R}f_{\mathbf{e}}(\xi,\tau')\widehat{\psi}(\tau-\tau')\cdot (\tau-\tau')\,d\tau'\Big|\Big|_{L^2_{\xi,\tau}}\leq C||f_{\mathbf{e}}||_{Y_k^{\mathbf{e}}}.$$
Thus the second term in \eqref{ni9} is bounded by $C||f_{\mathbf{e}}||_{Y_k^{\mathbf{e}}}$, which completes the proof of \eqref{ni8}.
\end{proof}

We prove now several nonlinear estimates. For $u\in C(\mathbb{R}:H^\infty)$ we define
\begin{equation}\label{no1}
\mathcal{N}(u)=\psi(t)\cdot 2\overline{u}(1+u\overline{u})^{-1}\sum_{j=1}^d(\partial_{x_j}u)^2\in C(\mathbb{R}:H^\infty),
\end{equation}
which is the nonlinear term in \eqref{Sch5}. We are looking to control
\begin{equation*}
\|\mathcal{N}(u)-\mathcal{N}(v)\|_{N^\sigma},\,\sigma>(d+1)/2,
\end{equation*}
where $u,v\in F^\sigma$. The plan is the following: if we ignore the factor $2(1+u\overline{u})^{-1}$, then $\mathcal{N}(u)$ is essentially of the form
\begin{equation*}
\psi(t)\cdot\overline{u}_1\cdot \nabla_xu_2\cdot\nabla_xu_3.
\end{equation*}
This is a trilinear expression. To estimate it, we use Lemma \ref{Lemmaa1} and the restriction $\sigma>(d+1)/2$ to place the two low-frequency factors in $L^{2,\infty}_\mathbf{e}$, for suitable vectors $\mathbf{e}$. Then, using Lemma \ref{Lemmas2}, we place the high frequency factor in $L^{\infty,2}_{\mathbf{e}}$, and gain $1/2$ derivative. The product is then in $L^{1,2}_{\mathbf{e}}$, which gains the second $1/2$ derivative (compare with the definition \eqref{v2}). 

There are certain technical difficulties to running this argument, mostly due to the presence of the factor $(1+u\overline{u})^{-1}$ and the fact that the spaces $F^\sigma$ are not stable under complex conjugation. To address this last problem (see \eqref{no50} below), we define normed spaces $\widetilde{Z}_k$, $k\in\mathbb{Z}_+$, and $\widetilde{F}^\sigma$, $\sigma\geq 0$:
\begin{equation}\label{no2}
\begin{split}
\widetilde{Z}_k&=\{f\in L^2(\mathbb{R}^d\times\mathbb{R}):f\text{ supported in }D_{k,\infty}\text{ and }\|f\|_{\widetilde{Z}_k}<\infty\},
\end{split}
\end{equation}  
where
\begin{equation}\label{no3}
\begin{split}
\|f\|_{\widetilde{Z}_k}&=2^{k/2}\sup_{\mathbf{e}\in\mathbb{S}^{d-1}}\|\mathcal{F}_{(d+1)}^{-1}[f\cdot \chi_{k,20}(\xi\cdot\mathbf{e})]\|_{L^{\infty,2}_{\mathbf{e}}}\\
&+2^{-(d-1)k/2}(k+1)^{-2}\sup_{\mathbf{e}\in\mathbb{S}^{d-1}}\|\mathbf{1}_{[-2,2]}(t)\cdot \mathcal{F}_{(d+1)}^{-1}(f)\|_{L^{2,\infty}_{\mathbf{e}}},
\end{split}
\end{equation}
and
\begin{equation}\label{no4}
\widetilde{F}^\sigma=\{u\in C(\mathbb{R}:H^\infty):\|u\|_{\widetilde{F}^\sigma}^2=\sum_{k=0}^\infty 2^{2\sigma k}\|\eta_k^{(d)}(\xi)\cdot \mathcal{F}_{(d+1)}u\|_{\widetilde{Z}_k}^2<\infty\}.
\end{equation}
In view of Lemma \ref{Lemmas2} and Lemma \ref{Lemmaa1},
\begin{equation}\label{no8}
\begin{cases}
&\|f\|_{\widetilde{Z}_k}\leq C\|f\|_{Z_k}\text{ for any }k\in\mathbb{Z}_+\text{ and }f\in Z_k;\\
&\|u\|_{\widetilde{F}^\sigma}\leq C\|u\|_{F^\sigma}\text{ for any }\sigma\geq 0\text{ and }u\in F^\sigma.
\end{cases}
\end{equation}
In addition, directly from the definition,
\begin{equation}\label{no50}
\|\overline{u}\|_{\widetilde{F}^\sigma}=\|u\|_{\widetilde{F}^\sigma}\text{ for any }\sigma\geq 0\text{ and }u\in \widetilde{F}^\sigma. 
\end{equation}
We start with a symmetric trilinear estimate. For $\sigma\in\mathbb{R}$ let $J^\sigma$ denote the operator defined by the Fourier multiplier $(\xi,\tau)\to(1+|\xi|^2)^{\sigma/2}$.

\newtheorem{Lemmaq4}[Lemmaq1]{Lemma}
\begin{Lemmaq4}\label{Lemmaq4}
If $\sigma>(d+1)/2$ and $u_1,u_2,u_3\in \widetilde{F}^\sigma$ then $\psi(t)\cdot J^1(u_1)\cdot J^1(u_2)\cdot J^1(u_3)\in N^\sigma$ and 
\begin{equation}\label{no11}
\begin{split}
\|\psi(t)\cdot J^1(u_1)\cdot J^1(u_2)\cdot J^1(u_3)\|_{N^{\sigma}}\leq C_\sigma\|u_1\|_{\widetilde{F}^\sigma}\cdot \|u_2\|_{\widetilde{F}^\sigma}\cdot \|u_3\|_{\widetilde{F}^\sigma}.
\end{split}
\end{equation}
\end{Lemmaq4}

\begin{proof}[Proof of Lemma \ref{Lemmaq4}] We fix a smooth function $\gamma^{(d)}:\mathbb{R}^d\to[0,1]$ supported in $[-2/3,2/3]^d$, equal to $1$ in $[-1/3,1/3]^d$, with the property that
\begin{equation*}
\sum_{m\in\mathbb{Z}^d}\gamma^{(d)}(\xi-m)\equiv 1.
\end{equation*}
Let $U=\psi(t)\cdot J^1(u_1)\cdot J^1(u_2)\cdot J^1(u_3)$. Using the definitions, 
\begin{equation}\label{no13}
\begin{split}
\|U\|_{N^{\sigma}}^2&\leq  C\sum_{k\in\mathbb{Z}_+}2^{2\sigma k}\sum_{|m|\in[2^{8d},2^{12d}]}\\
&\|\eta_k^{(d)}(\xi)\cdot \gamma^{(d)}(\xi/2^{k-10d}-m)\cdot (\tau+|\xi|^2+i)^{-1}\cdot\mathcal{F}_{(d+1)}(U)\|_{Z_k}^2.
\end{split}
\end{equation}
For $k\in\mathbb{Z}_+$ let $Q_k$ denote the operator defined by the Fourier multiplier $(\xi,\tau)\to\eta_k^{(d)}(\xi)$. We have
\begin{equation}\label{no15}
Q_k[Q_{k_1}(v_1)\cdot Q_{k_2}(v_2)\cdot Q_{k_3}(v_3)]\equiv 0\text{ unless }k\leq \max(k_1,k_2,k_3)+3.
\end{equation}
In view of \eqref{no13}, for \eqref{no11} it suffices to prove that for $|m|\in[2^{8d},2^{12d}]$ fixed
\begin{equation}\label{no14}
\begin{split}
\sum_{k\in\mathbb{Z}_+}2^{2\sigma k}\|\gamma^{(d)}(\xi/2^{k-10d}-m)\cdot (\tau+|\xi|^2&+i)^{-1}\cdot\mathcal{F}_{(d+1)}(Q_k(U))\|_{Z_k}^2\\
&\leq C_\sigma\|u_1\|_{\widetilde{F}^\sigma}^2\cdot \|u_2\|_{\widetilde{F}^\sigma}^2\cdot \|u_3\|_{\widetilde{F}^\sigma}^2.
\end{split}
\end{equation}
Let $\widehat{m}=m/|m|$ and define
\begin{equation*}
S(\widehat{m})=\{\mathbf{e}\in\{\mathbf{e}_1,\ldots,\mathbf{e}_L\}:|\mathbf{e}\cdot\widehat{m}|\geq 3/4\}\text{ and }\widetilde{L}_{\widehat{m}}^{p,q}=\oplus_{\mathbf{e}\in S(\widehat{m})}L^{p,q}_{\mathbf{e}}.
\end{equation*}
Using the definition of $Z_k$  (and \eqref{gu76} if $k\leq 99$), and the identity $L^{p,q}_\mathbf{e}=L^{p,q}_{-\mathbf{e}}$,
\begin{equation}\label{no20}
\begin{split}
\|\gamma^{(d)}(\xi/2^{k-10d}-m)\cdot &(\tau+|\xi|^2+i)^{-1}\cdot\mathcal{F}_{(d+1)}(Q_k(U))\|_{Z_k}\\
&\leq C2^{-k/2}\|Q_k(U)\|_{\widetilde{L}^{1,2}_{\widehat{m}}}.
\end{split}
\end{equation}

We assume now that $k$ is fixed and estimate the right-hand side of \eqref{no20}. In view of \eqref{no15} and the definition of $U$,
\begin{equation*}
Q_k(U)=Q_k\Big[\psi(t) \sum_{(k_1,k_2,k_3)\in T_k}J^1Q_{k_1}(u_1)\cdot J^1Q_{k_2}(u_2)\cdot J^1Q_{k_3}(u_3) \Big],
\end{equation*}
where $T_k=\{(k_1,k_2,k_3)\in(\mathbb{Z}_+)^3:k\leq \max(k_1,k_2,k_3)+3\}$. Since $Q_k$ is a bounded operator on $\widetilde{L}_{\widehat{m}}^{p,q}$ uniformly in $k$, the right-hand side of \eqref{no20} is dominated by
\begin{equation}\label{no39}
C2^{-k/2}\sum_{(k_1,k_2,k_3)\in T_k}\big|\big|\psi(t)\cdot J^1Q_{k_1}(u_1)\cdot J^1Q_{k_2}(u_2)\cdot J^1Q_{k_3}(u_3) \big|\big|_{\widetilde{L}^{1,2}_{\widehat{m}}}.
\end{equation}

Assume, by symmetry, that $k_1=\max(k_1,k_2,k_3)\geq k-4$. Then, using the definitions \eqref{no3}  and \eqref{no4}, for any vector $\mathbf{e}\in\mathbb{S}^{d-1}$,
\begin{equation}\label{no40}
\begin{split}
\|\mathbf{1}_{[-2,2]}(t)\cdot J^1Q_{k_2}(u_2)\|_{L^{2,\infty}_\mathbf{e}}&\leq C2^{k_2}\|\mathbf{1}_{[-2,2]}(t)\cdot Q_{k_2}(u_2)\|_{L^{2,\infty}_\mathbf{e}}\\
&\leq C2^{-\sigma k_2}2^{(d+1)k_2/2}(k_2+1)^2\|u_2\|_{\widetilde{F}^\sigma}.
\end{split}
\end{equation}
Similarly, for any vector $\mathbf{e}\in\mathbb{S}^{d-1}$,
\begin{equation}\label{no41}
\begin{split}
\|\mathbf{1}_{[-2,2]}(t)\cdot J^1Q_{k_3}(u_3)\|_{L^{2,\infty}_\mathbf{e}}\leq C2^{-\sigma k_3}2^{(d+1)k_3/2}(k_3+1)^2\|u_3\|_{\widetilde{F}^\sigma}.
\end{split}
\end{equation}
We show next that
\begin{equation}\label{no42}
\|J^1Q_{k_1}(u_1)\|_{\widetilde{L}^{\infty,2}_{\widehat{m}}}\leq C2^{k_1/2}||\mathcal{F}_{(d+1)}(Q_{k_1}(u_1))||_{\widetilde{Z}_{k_1}}.
\end{equation}
Using the definition \eqref{no3}, for \eqref{no42} it suffices to prove that
\begin{equation}\label{no43}
\|Q_{k_1}(u_1)\|_{\widetilde{L}^{\infty,2}_{\widehat{m}}}\leq C\sup_{\mathbf{e}\in\mathbb{S}^{d-1}}||\mathcal{F}^{-1}_{(d+1)}[\mathcal{F}_{(d+1)}(Q_{k_1}(u_1))\cdot \eta_{k_1,20}(\xi\cdot\mathbf{e})]||_{L^{\infty,2}_\mathbf{e}}.
\end{equation}
Since $\chi_{k_1,20}\equiv  1$ if $k_1\leq 99$, we may assume $k_1\geq 100$ in \eqref{no43}. Using the function $\gamma^{(d)}$ defined at the beginning of the proof, we decompose
\begin{equation*}
Q_{k_1}(u_1)=\sum_{|n|\in[2^{8d},2^{12d}]}\mathcal{F}^{-1}_{(d+1)}[\mathcal{F}_{(d+1)}(Q_{k_1}(u_1))\cdot\gamma^{(d)}(\xi/2^{k_1-10d}-n)].
\end{equation*}
Thus, for \eqref{no43} we only need the following elementary statement: if  $\widehat{n}\in\mathbb{S}^{d-1}$ there is $\mathbf{e}\in S(\widehat{m})$ with the property that $\mathbf{e}\cdot\widehat{n}\geq 2^{-10}$. To see this, we find first a vector $\mathbf{e}'\in\mathbb{S}^{d-1}$ with the properties $|\mathbf{e}'\cdot\widehat{m}|\geq 7/8$ and $\mathbf{e}'\cdot\widehat{n}\geq 2^{-9}$ (simply take $\mathbf{e'}=\widehat{m}$ or $\mathbf{e}'=-\widehat{m}$ or $\mathbf{e}'=\widehat{\widehat{m}+2^{-6}\widehat{n}}$), and  then find a vector $\mathbf{e}\in\{\mathbf{e}_1,\ldots,\mathbf{e}_L\}$ such that $|\mathbf{e}-\mathbf{e}'|\leq 2^{-50}$. This completes the proof of \eqref{no42}.

Using \eqref{no40}, \eqref{no41}, \eqref{no42}, and the restriction $\sigma>(d+1)/2$, and summing over $k_2,k_3$, the expression in \eqref{no39} is bounded by
\begin{equation*}
\begin{split}
&C_\sigma||u_2||_{\widetilde{F}^\sigma}\cdot ||u_3||_{\widetilde{F}^\sigma}\sum_{k_1\geq k-4}2^{(k_1-k)/2}||\mathcal{F}_{(d+1)}(Q_{k_1}(u_1))||_{\widetilde{Z}_{k_1}}\\
+&C_\sigma||u_1||_{\widetilde{F}^\sigma}\cdot ||u_3||_{\widetilde{F}^\sigma}\sum_{k_2\geq k-4}2^{(k_2-k)/2}||\mathcal{F}_{(d+1)}(Q_{k_2}(u_2))||_{\widetilde{Z}_{k_2}}\\
+&C_\sigma||u_1||_{\widetilde{F}^\sigma}\cdot ||u_2||_{\widetilde{F}^\sigma}\sum_{k_3\geq k-4}2^{(k_3-k)/2}||\mathcal{F}_{(d+1)}(Q_{k_3}(u_3))||_{\widetilde{Z}_{k_3}}.
\end{split}
\end{equation*}
The bound  \eqref{no14} follows from \eqref{no20}, which completes the proof of  Lemma \ref{Lemmaq4}.
\end{proof}

We continue with a symmetric multilinear estimate.

\newtheorem{Lemmaq5}[Lemmaq1]{Lemma}
\begin{Lemmaq5}\label{Lemmaq5}
If $\sigma>(d+1)/2$, $n\geq 1$, and $u_1,\ldots,u_n\in F^\sigma$ then $\widetilde{u}_1\cdot\ldots\cdot \widetilde{u}_n\in\widetilde{F}^\sigma$ and 
\begin{equation}\label{no60}
||\widetilde{u}_1\cdot\ldots\cdot \widetilde{u}_n||_{\widetilde{F}^\sigma}\leq (C_\sigma)^n\cdot ||u_1||_{F^\sigma}\cdot\ldots\cdot ||u_n||_{F^\sigma},
\end{equation}
where $\widetilde{u}_m\in\{u_m,\overline{u}_m\}$ for $m=1,\ldots,n$.
\end{Lemmaq5}

\begin{proof}[Proof of Lemma \ref{Lemmaq5}] In view  of  \eqref{no8} and  \eqref{no50}, we may assume $n\geq 2$. We recall the definition
\begin{equation}\label{no62}
||\widetilde{u}_1\cdot\ldots\cdot \widetilde{u}_n||_{\widetilde{F}^\sigma}^2=\sum_{k\in\mathbb{Z}_+}2^{2\sigma k}||\mathcal{F}_{(d+1)}[Q_k(\widetilde{u}_1\cdot\ldots\cdot \widetilde{u}_n)]||_{\widetilde{Z}_k}^2.
\end{equation}
We have 
\begin{equation*}
Q_k[Q_{k_1}(\widetilde{u}_1)\cdot\ldots\cdot Q_{k_n}(\widetilde{u}_n)]=0\text{ unless }\max(k_1,\ldots,k_n)\geq k-2-\log_2n.
\end{equation*}
We fix  $k\in\mathbb{Z}_+$ and  let $T_k^n=\{(k_1,\ldots,k_n)\in(\mathbb{Z}_+)^n:k\leq \max(k_1,\ldots,k_n)+2+\log n\}$. Then
\begin{equation}\label{no63}
||\mathcal{F}_{(d+1)}[Q_k(\widetilde{u}_1\cdot\ldots\cdot \widetilde{u}_n)]||_{\widetilde{Z}_k}\leq \sum_{(k_1,\ldots,k_n)\in T_k^n}||\mathcal{F}_{(d+1)}[Q_k(Q_{k_1}(\widetilde{u}_1)\cdot\ldots\cdot Q_{k_n}(\widetilde{u}_n))]||_{\widetilde{Z}_k}.
\end{equation}

To analyze the right-hand side of \eqref{no63} for $(k_1,\ldots,k_n)$ fixed, assume, by symmetry, that $k_1=\max(k_1,\ldots,k_n)$ and $k_2=\max(k_2,\ldots,k_n)$. Using \eqref{lb44}, Lemma \ref{Lemmaa1},  the fact that $\sigma>(d+1)/2$, and examining the definition \eqref{no3}
\begin{equation}\label{no64}
\begin{split}
&2^{-(d-1)k/2}(k+1)^{-2}\sup_{\mathbf{e}\in\mathbb{S}^{d-1}}||\mathbf{1}_{[-2,2]}(t)\cdot Q_k[Q_{k_1}(\widetilde{u}_1)\cdot\ldots\cdot Q_{k_n}(\widetilde{u}_n)]||_{L^{2,\infty}_{\mathbf{e}}}\\
&\leq C2^{-(d-1)k/2}(k+1)^{-2}\sup_{\mathbf{e}\in\mathbb{S}^{d-1}}||\mathbf{1}_{[-2,2]}(t)Q_{k_1}(\widetilde{u}_1)||_{L^{2,\infty}_{\mathbf{e}}}\cdot\prod_{m=2}^n||Q_{k_m}(\widetilde{u}_m)||_{L^\infty}\\
&\leq C^n\frac{2^{(d-1)k_1/2}(k_1+1)^2}{2^{(d-1)k/2}(k+1)^2}||\mathcal{F}_{(d+1)}(Q_{k_1}(u_1))||_{Z_{k_1}}\cdot \prod_{m=2}^n2^{dk_m/2}||\mathcal{F}_{(d+1)}(Q_{k_m}(u_m))||_{Z_{k_m}}.
\end{split}
\end{equation}
To estimate the $L^{\infty,2}_{\mathbf{e}}$ norm in the first line of \eqref{no3} we consider two cases. If $k_2\geq k_1-50\log_2n$ then, using  \eqref{lb44}, \eqref{gu90}, and the restriction $\sigma>(d+1)/2$
\begin{equation}\label{no65}
\begin{split}
&2^{k/2}\sup_{\mathbf{e}\in\mathbb{S}^{d-1}}||\mathcal{F}_{(d+1)}^{-1}\{\mathcal{F}_{(d+1)}[Q_k(Q_{k_1}(\widetilde{u}_1)\cdot\ldots\cdot Q_{k_n}(\widetilde{u}_n))](\xi,\tau)\cdot \chi_{k,20}(\xi\cdot\mathbf{e})\}||_{L^{\infty,2}_{\mathbf{e}}}\\
&\leq C2^{k/2}\sup_{\mathbf{e}\in\mathbb{S}^{d-1}}||Q_{k_1}(\widetilde{u}_1)||_{L^{\infty,2}_{\mathbf{e}}}\cdot\prod_{m=2}^n||Q_{k_m}(\widetilde{u}_m)||_{L^\infty}\leq C^n||\mathcal{F}_{(d+1)}(Q_{k_1}(u_1))||_{Z_{k_1}}\\
&\times2^{(d+1)k_2/2}(k_2+1)||\mathcal{F}_{(d+1)}(Q_{k_2}(u_2))||_{Z_{k_2}}\cdot \prod_{m=3}^n2^{dk_m/2}||\mathcal{F}_{(d+1)}(Q_{k_m}(u_m))||_{Z_{k_m}}.
\end{split}
\end{equation}
If $k_2\leq k_1-50\log_2n$ (so $|k-k_1|\leq 2$), let $U=Q_{k_2}(\widetilde{u}_2)\cdot\ldots\cdot Q_{k_n}(\widetilde{u}_n)$ and notice that $\mathcal{F}_{(d+1)}(U)$ is supported in the set $\{(\xi,\tau):|\xi|\leq 2^{k_1-40}\}$. Thus, for any $\mathbf{e}\in\mathbb{S}^{d-1}$, 
\begin{equation*}
\mathcal{F}_{(d+1)}[Q_{k_1}(\widetilde{u}_1)\cdot U](\xi,\tau)\cdot\chi_{k,20}(\xi\cdot\mathbf{e})=\mathcal{F}_{(d+1)}[Q^\mathbf{e}_{k_1}(\widetilde{u}_1)\cdot U](\xi,\tau)\cdot\chi_{k,20}(\xi\cdot\mathbf{e})
\end{equation*}
where
\begin{equation*}
\mathcal{F}_{(d+1)}(Q_{k_1}^\mathbf{e}(u))(\xi,\tau)=\mathcal{F}_{(d+1)}(Q_{k_1}(u))(\xi,\tau)\cdot \chi_{k,30}(\xi\cdot\mathbf{e}).
\end{equation*}
Thus, using Lemma \ref{Lemmas2} and \eqref{lb44},
\begin{equation}\label{no66}
\begin{split}
&2^{k/2}\sup_{\mathbf{e}\in\mathbb{S}^{d-1}}||\mathcal{F}_{(d+1)}^{-1}\{\mathcal{F}_{(d+1)}[Q_k(Q_{k_1}(\widetilde{u}_1)\cdot\ldots\cdot Q_{k_n}(\widetilde{u}_n))](\xi,\tau)\cdot \chi_{k,20}(\xi\cdot\mathbf{e})\}||_{L^{\infty,2}_{\mathbf{e}}}\\
&\leq C2^{k/2}\sup_{\mathbf{e}\in\mathbb{S}^{d-1}}||Q_{k_1}^\mathbf{e}(\widetilde{u}_1)||_{L^{\infty,2}_{\mathbf{e}}}\cdot\prod_{m=2}^n||Q_{k_m}(\widetilde{u}_m)||_{L^\infty}\\
&\leq C^n||\mathcal{F}_{(d+1)}(Q_{k_1}(u_1))||_{Z_{k_1}}\cdot \prod_{m=2}^n2^{dk_m/2}||\mathcal{F}_{(d+1)}(Q_{k_m}(u_m))||_{Z_{k_m}}.
\end{split}
\end{equation}

We combine  \eqref{no64}, \eqref{no65}, and \eqref{no66}, and sum over $k_2,\ldots,k_n\in\mathbb{Z}_+$. It follows that the part of the expression in the right-hand side of \eqref{no63} which corresponds to $k_1=\max(k_1,\ldots,k_n)$ is dominated by
\begin{equation*}
(C_\sigma)^n\prod_{m=2}^n||u_m||_{F^\sigma}\cdot\sum_{k_1\geq k-2-\log_2n}\frac{2^{(d-1)k_1/2}(k_1+1)^2}{2^{(d-1)k/2}(k+1)^2}\cdot ||\mathcal{F}_{(d+1)}(Q_{k_1}(u_1))||_{Z_{k_1}}
\end{equation*}
The bound  \eqref{no60} then follows from \eqref{no62}.
\end{proof}

For $u\in C(\mathbb{R}:H^\infty)$ we define
\begin{equation}\label{no70}
\mathcal{N}_0(u)=2\overline{u}(1+u\overline{u})^{-1}\in C(\mathbb{R}:H^\infty),
\end{equation}
so $\mathcal{N}(u)=\psi(t)\cdot \mathcal{N}_0(u)\cdot \sum_{j=1}^d(\partial_{x_j}u)^2$ (compare with \eqref{no1}).

\newtheorem{Lemmaq6}[Lemmaq1]{Lemma}
\begin{Lemmaq6}\label{Lemmaq6}
Assume $\sigma>(d+1)/2$. Then there is $\overline{c}(\sigma)>0$ with the property that 
\begin{equation}\label{no71}
\|J^{\sigma'}(\mathcal{N}_0(u)-\mathcal{N}_0(v))\|_{\widetilde{F}^\sigma}\leq C(\sigma,\sigma',\|J^{\sigma'}u\|_{F^\sigma}+\|J^{\sigma'}v\|_{F^\sigma})\cdot \|J^{\sigma'}(u-v)\|_{F^\sigma} 
\end{equation}
for any $\sigma'\in\mathbb{Z}_+$, and any $u,v\in B_{F^\sigma}(0,\overline{c}(\sigma))\cap F^{\sigma+\sigma'}$.
\end{Lemmaq6}

\begin{proof}[Proof of Lemma \ref{Lemmaq6}] We write first
\begin{equation*}
\mathcal{N}_0(u)-\mathcal{N}_0(v)=\frac{\overline{u-v}}{(1+u\overline{u})(1+v\overline{v})}-\frac{(u-v)\cdot \overline{u}\,\overline{v}}{(1+u\overline{u})(1+v\overline{v})}.
\end{equation*}
By expanding in power series and using \eqref{no50}, it suffices to prove that
\begin{equation}\label{no72}
\|J^{\sigma'}((u-v)\cdot u_1\cdot \ldots\cdot u_n)\|_{\widetilde{F}^\sigma}\leq 2^{-n}\cdot C(\sigma,\sigma',\|J^{\sigma'}u\|_{F^\sigma}+\|J^{\sigma'}v\|_{F^\sigma})\cdot \|J^{\sigma'}(u-v)\|_{F^\sigma}, 
\end{equation}
for any $n\geq 1$, where $u_m\in\{u,\overline{u},v,\overline{v}\}$. This follows directly from Lemma \ref{Lemmaq5}: since $\sigma'\in\mathbb{Z}_+$ we can distribute the $\sigma'$ derivatives in the left-hand side of \eqref{no72} in at most $(n+1)^{\sigma'}\leq 2^n\cdot C_{\sigma'}$ ways. For each of the resulting terms we use \eqref{no60}; since at most $\sigma'$ of the factors contain derivatives, all the other $n-\sigma'$ factors contribute a factor of $\overline{c}(\sigma)\ll 1$, which gives the exponential decay in \eqref{no72}.
\end{proof}

\section{Proof of Theorem \ref{Main2}}\label{proof}

In this section we complete the proof of Theorem \ref{Main2}. Our  main ingredients are Lemma \ref{Lemmaq1}, Lemma \ref{Lemmaq3}, Lemma \ref{Lemmaq4}, Lemma \ref{Lemmaq6}, and the bound \begin{equation}\label{fi1}
\sup_{t\in\mathbb{R}}\|u\|_{H^\sigma}\leq C_\sigma\|u\|_{F^\sigma}\text{ for any }\sigma\geq 0\text{ and }u\in F^\sigma,
\end{equation}
which follows from Lemma \ref{Lemmaa2}. Assume, as in Theorem \ref{Main2} that ${\sigma_0}>(d+1)/2$ and $\phi\in H^\infty\cap B_{H^{\sigma_0}}(0,\epsilon({\sigma_0}))$, where $\epsilon({\sigma_0})\ll 1$ is to be fixed. We define recursively
\begin{equation}\label{fi2}
\begin{cases}
&u_0=\psi(t)\cdot W(t)\phi;\\
&u_{n+1}=\psi(t)\cdot W(t)\phi+\psi(t)\cdot \int_0^tW(t-s)(\mathcal{N}(u_n(s)))\,ds\text{ for }n\in\mathbb{Z}_+.
\end{cases}
\end{equation}
Clearly, $u_n\in C(\mathbb{R}:H^\infty)$. 

We show first that
\begin{equation}\label{fi3}
\|u_n\|_{F^{\sigma_0}}\leq C_{\sigma_0}\|\phi\|_{H^{\sigma_0}}\text{ for any }n=0,1,\ldots,\text{ if }\epsilon({\sigma_0}) \text{ is sufficiently small}. 
\end{equation} 
The bound \eqref{fi3} holds for $n=0$, due to Lemma \ref{Lemmaq1}. Then, using Lemma \ref{Lemmaq6} with $\sigma'=0$, $v\equiv 0$, Lemma \ref{Lemmaq4}, and the inequality \eqref{no8}, we have
\begin{equation*}
\|\mathcal{N}(u_n)\|_{N^{\sigma_0}}\leq C_{\sigma_0}\|u_n\|_{F^{\sigma_0}}^3.
\end{equation*}
Using Lemma \ref{Lemmaq3}, the definition \eqref{fi2}, and Lemma \ref{Lemmaq1}, it follows that
\begin{equation*}
\|u_{n+1}\|_{F^{\sigma_0}}\leq C_{\sigma_0}\|\phi\|_{H^{\sigma_0}}+C_{\sigma_0}\|u_n\|_{F^{\sigma_0}}^3,
\end{equation*}
which leads to \eqref{fi3} by induction over $n$.

We show now that 
\begin{equation}\label{fi4}
\|u_{n}-u_{n-1}\|_{F^{\sigma_0}}\leq 2^{-n}\cdot C_{\sigma_0}\|\phi\|_{H^{\sigma_0}}\text{ for any }n\in\mathbb{Z}_+\text{ if }\epsilon({\sigma_0}) \text{ is sufficiently small}. 
\end{equation}
This is clear for $n=0$ (with $u_{-1}\equiv 0$), using Lemma \ref{Lemmaq1}. Then, using Lemma \ref{Lemmaq6} with $\sigma'=0$, Lemma \ref{Lemmaq4}, and the estimates \eqref{no8} and \eqref{fi3}, we have
\begin{equation*}
\|\mathcal{N}(u_{n-1})-\mathcal{N}(u_{n-2})\|_{N^{\sigma_0}}\leq C_{\sigma_0}\cdot \epsilon({\sigma_0})^2\cdot \|u_{n-1}- u_{n-2}\|_{F^{\sigma_0}}.
\end{equation*}
Using Lemma \ref{Lemmaq3} and the definition \eqref{fi2} it follows that
\begin{equation*}
\|u_{n}-u_{n-1}\|_{F^{\sigma_0}}\leq C_{\sigma_0}\cdot \epsilon({\sigma_0})^2\cdot \|u_{n-1}-u_{n-2}\|_{F^{\sigma_0}},
\end{equation*}
which leads to \eqref{fi4} by induction over $n$.

We show now that 
\begin{equation}\label{fi5}
\|J^{\sigma'}(u_n)\|_{F^{\sigma_0}}\leq C({\sigma_0},\sigma',\|J^{\sigma'}\phi\|_{H^{\sigma_0}})\text{ for any }n,\sigma'\in\mathbb{Z}_+.
\end{equation}
We argue by induction over $\sigma'$ (the case $\sigma'=0$ follows from \eqref{fi3}). So we may assume that
\begin{equation}\label{fi6}
\|J^{\sigma'-1}(u_n)\|_{F^{\sigma_0}}\leq C({\sigma_0},\sigma',\|J^{\sigma'-1}\phi\|_{H^{\sigma_0}})\text{ for any }n\in\mathbb{Z}_+,
\end{equation}
and it suffices to prove that
\begin{equation}\label{fi7}
\|\partial_{x_i}^{\sigma'}(u_n)\|_{F^{\sigma_0}}\leq C({\sigma_0},\sigma',\|J^{\sigma'}\phi\|_{H^{\sigma_0}})\text{ for any }n\in\mathbb{Z}_+\text{ and }i=1,\ldots,d.
\end{equation}
The bound \eqref{fi7} for $n=0$ follows from Lemma \ref{Lemmaq1}. We use the decomposition
\begin{equation*}
\mathcal{N}(u_n)=\sum_{j=1}^d\psi(t)\cdot\mathcal{N}_0(u_n)\cdot(\partial_{x_j}u_n)^2,
\end{equation*}
thus
\begin{equation}\label{fi8}
\partial_{x_i}^{\sigma'}(\mathcal{N}(u_n))=2\sum_{j=1}^d\psi(t)\cdot\mathcal{N}_0(u_n)\cdot\partial_{x_j}u_n\cdot\partial_{x_i}^{\sigma'}\partial_{x_j}u_n+E_n,
\end{equation}
where
\begin{equation*}
E_n=\sum_{j=1}^d\psi(t)\cdot \sum_{\sigma'_1+\sigma'_2+\sigma'_3=\sigma'\text{ and }\sigma'_3,\sigma'_2<\sigma'}\partial_{x_1}^{\sigma'_1}\mathcal{N}_0(u_n)\cdot\partial_{x_i}^{\sigma'_2}\partial_{x_j}u_n\cdot\partial_{x_i}^{\sigma'_3}\partial_{x_j}u_n.
\end{equation*}
Using  Lemma \ref{Lemmaq4}, 
\begin{equation*}
\begin{split}
||E_n||_{N^{\sigma_0}}\leq C_{\sigma_0}&\sum_{j=1}^d\sum_{\sigma'_1+\sigma'_2+\sigma'_3=\sigma'\text{ and }\sigma'_3,\sigma'_2<\sigma'}\\
&||J^{-1}\partial_{x_1}^{\sigma'_1}\mathcal{N}_0(u_n)||_{\widetilde{F}^{\sigma_0}}\cdot ||J^{-1}\partial_{x_i}^{\sigma'_2}\partial_{x_j}u_n||_{\widetilde{F}^{\sigma_0}}\cdot ||J^{-1}\partial_{x_i}^{\sigma'_3}\partial_{x_j}u_n||_{\widetilde{F}^{\sigma_0}}.
\end{split}
\end{equation*}
Using now Lemma \ref{Lemmaq6} with $v=0$, the bound \eqref{no8}, and the induction hypothesis \eqref{fi6}, we have
\begin{equation}\label{fi9}
||E_n||_{N^{\sigma_0}}\leq C({\sigma_0},\sigma',\|J^{\sigma'-1}\phi\|_{H^{\sigma_0}}).
\end{equation}
In addition, using again Lemma \ref{Lemmaq4}, Lemma \ref{Lemmaq6} with $v=0$, \eqref{no8} and \eqref{fi3},
\begin{equation}\label{fi10}
||2\sum_{j=1}^d\psi(t)\cdot\mathcal{N}_0(u_n)\cdot\partial_{x_j}u_n\cdot\partial_{x_i}^{\sigma'}\partial_{x_j}u_n||_{N^{\sigma_0}}\leq C_{{\sigma_0}}\cdot \epsilon({\sigma_0})^2\cdot ||\partial_{x_i}^{\sigma'}u_n||_{F^{\sigma_0}}.
\end{equation}
We use now the definition \eqref{fi2}, together with Lemma \ref{Lemmaq1}, Lemma \ref{Lemmaq3}, and the bounds \eqref{fi9} and \eqref{fi10} to conclude that
\begin{equation*}
||\partial_{x_i}^{\sigma'}u_{n+1}||_{F^{\sigma_0}}\leq C({\sigma_0},\sigma',\|J^{\sigma'}\phi\|_{H^{\sigma_0}})+C_{{\sigma_0}}\cdot \epsilon({\sigma_0})^2\cdot ||\partial_{x_i}^{\sigma'}u_n||_{F^{\sigma_0}}.
\end{equation*}
The bound  \eqref{fi7} follows by induction over $n$ provided that $\epsilon({\sigma_0})$ is sufficiently small.

Finally, we show that
\begin{equation}\label{fi20}
\|J^{\sigma'}(u_{n}-u_{n-1}))\|_{F^{\sigma_0}}\leq 2^{-n}\cdot C({\sigma_0},\sigma',\|J^{\sigma'}\phi\|_{H^{\sigma_0}})\text{ for any }n,\sigma'\in\mathbb{Z}_+.
\end{equation}
As before, we argue by induction over $\sigma'$ (the case $\sigma'=0$ follows from \eqref{fi4}). So we may assume that
\begin{equation}\label{fi61}
\|J^{\sigma'-1}(u_n-u_{n-1})\|_{F^{\sigma_0}}\leq 2^{-n}\cdot C({\sigma_0},\sigma',\|J^{\sigma'-1}\phi\|_{H^{\sigma_0}})\text{ for any }n\in\mathbb{Z}_+,
\end{equation}
and it suffices to prove that
\begin{equation}\label{fi71}
\|\partial_{x_i}^{\sigma'}(u_n-u_{n-1})\|_{F^{\sigma_0}}\leq 2^{-n}\cdot C({\sigma_0},\sigma',\|J^{\sigma'}\phi\|_{H^{\sigma_0}})\text{ for any }n\in\mathbb{Z}_+\text{ and }i=1,\ldots,d.
\end{equation}
The bound \eqref{fi71} for $n=0$ follows from Lemma \ref{Lemmaq1}. For $n\geq 1$ we use the decomposition
\begin{equation}\label{fi73}
\begin{split}
\mathcal{N}(u_{n-1})-&\mathcal{N}(u_{n-2})=\sum_{j=1}^d\psi(t)\cdot(\mathcal{N}_0(u_{n-1})-\mathcal{N}_0(u_{n-2}))\cdot(\partial_{x_j}u_{n-1})^2\\
&+\sum_{j=1}^d\psi(t)\cdot\mathcal{N}_0(u_{n-2})\cdot\partial_{x_j}(u_{n-1}-u_{n-2})\cdot \partial_{x_j}(u_{n-1}+u_{n-2}).
\end{split}
\end{equation}
The same argument as before, which consists of expanding the $\sigma'$ derivative, and combining Lemma \ref{Lemmaq4}, Lemma \ref{Lemmaq6}, \eqref{fi5}, and \eqref{fi61}, shows that
\begin{equation}\label{fi72}
\begin{split}
\big|\big|\partial_{x_1}^{\sigma'}\big[\sum_{j=1}^d\psi(t)\cdot(\mathcal{N}_0(u_{n-1})-\mathcal{N}_0(u_{n-2}))\cdot(\partial_{x_j}u_{n-1})^2\big]\big|\big|_{N^{\sigma_0}}\\
\leq 2^{-n}\cdot C({\sigma_0},\sigma',||J^{\sigma'}\phi||_{H^{\sigma_0}}).
\end{split}
\end{equation}
To estimate the $\sigma'$ derivative of the term in the second line of \eqref{fi73}, we expand again the $\sigma'$ derivatives. Using again the combination of Lemma \ref{Lemmaq4}, Lemma \ref{Lemmaq6}, \eqref{fi5}, and \eqref{fi61}, the $N^{\sigma_0}$ norm of most of the terms that appear is again dominated by $2^{-n}\cdot C({\sigma_0},\sigma',||J^{\sigma'}\phi||_{H^{\sigma_0}})$. The only remaining terms are
\begin{equation*}
\sum_{j=1}^d\psi(t)\cdot\mathcal{N}_0(u_{n-2})\cdot\partial_{x_1}^{\sigma'}\partial_{x_j}(u_{n-1}-u_{n-2})\cdot \partial_{x_j}(u_{n-1}+u_{n-2}),
\end{equation*}
and we can estimate
\begin{equation*}
\begin{split}
||\sum_{j=1}^d\psi(t)\cdot\mathcal{N}_0(u_{n-2})\cdot\partial_{x_1}^{\sigma'}\partial_{x_j}(u_{n-1}-u_{n-2})\cdot \partial_{x_j}(u_{n-1}+u_{n-2})||_{N^{\sigma_0}}\\
\leq C_{\sigma_0}\cdot\epsilon({\sigma_0})^2\cdot ||\partial_{x_1}^{\sigma'}(u_{n-1}-u_{n-2})||_{F^{\sigma_0}}.
\end{split}
\end{equation*}
As before, it follows that
\begin{equation*}
\begin{split}
||\partial_{x_i}^{\sigma'}(u_{n}-u_{n-1})||_{F^{\sigma_0}}&\leq 2^{-n}\cdot C({\sigma_0},\sigma',\|J^{\sigma'}\phi\|_{H^{\sigma_0}})\\
&+C_{{\sigma_0}}\cdot \epsilon({\sigma_0})^2\cdot ||\partial_{x_i}^{\sigma'}(u_{n-1}-u_{n-2})||_{F^{\sigma_0}}.
\end{split}
\end{equation*}
The bound  \eqref{fi71} follows by induction provided that $\epsilon({\sigma_0})$ is sufficiently small.

We can now use \eqref{fi20} and \eqref{fi1} to construct
\begin{equation*}
u=\lim_{n\to\infty}u_n\in C(\mathbb{R}:H^\infty).
\end{equation*}
In view of \eqref{fi2}, 
\begin{equation*}
u=\psi(t)\cdot W(t)\phi+\psi(t)\cdot \int_0^tW(t-s)(\mathcal{N}(u(s)))\,ds\text{ on }\mathbb{R}^d\times\mathbb{R},
\end{equation*}
so $\widetilde{S}^\infty(\phi)$, the restriction of $u$ to $\mathbb{R}^d\times[-1,1]$, is a solution of the initial-value problem \eqref{Sch5}. The bound \eqref{by1} follows from the uniform bound \eqref{fi5} and  \eqref{fi1}.

For Theorem \ref{Main2} (b) and (c), it suffices to show that if $\sigma'\in\mathbb{Z}_+$ and $\phi,\phi'\in B_{H^{\sigma_0}}(0,\epsilon({\sigma_0}))\cap H^\infty$ then
\begin{equation}\label{fi80}
\sup_{t\in[-1,1]}||\widetilde{S}^\infty(\phi)-\widetilde{S}^\infty(\phi')||_{H^{\sigma_0+\sigma'}}\leq C(\sigma_0,\sigma',||\phi||_{H^{\sigma_0+\sigma'}})\cdot ||\phi-\phi'||_{H^{\sigma_0+\sigma'}}.
\end{equation}
Part  (b) corresponds to the case $\sigma'=0$. To prove \eqref{fi80}, we define the sequences $u_n$ and $u'_n$, $n\in\mathbb{Z}_+$, as in \eqref{fi2}. Using  Lemma \ref{Lemmaq1},
\begin{equation*}
||u_0-u'_0||_{F^{\sigma_0}}\leq C_{\sigma_0}||\phi-\phi'||_{H^{\sigma_0}}.
\end{equation*}
Then we decompose $\mathcal{N}(u_n)-\mathcal{N}(u'_n)$ as in \eqref{fi73}. As before, we combine Lemma \ref{Lemmaq1}, Lemma \ref{Lemmaq3}, Lemma \ref{Lemmaq4}, Lemma \ref{Lemmaq6}, and the uniform bound \eqref{fi3} to conclude that
\begin{equation*}
||u_{n+1}-u'_{n+1}||_{F^{\sigma_0}}\leq C_{\sigma_0}||\phi-\phi'||_{H^{\sigma_0}}+C_{\sigma_0}\cdot\epsilon({\sigma_0})^2\cdot ||u_{n}-u'_{n}||_{F^{\sigma_0}}.
\end{equation*}
By induction over $n$ it follows that
\begin{equation*}
||u_n-u'_n||_{F^{\sigma_0}}\leq C_{\sigma_0}||\phi-\phi'||_{H^{\sigma_0}}\text{ for any }n\in\mathbb{Z}_+.
\end{equation*}
In view of \eqref{fi1} this proves \eqref{fi80} for  $\sigma'=0$.

Assume now that $\sigma'\geq  1$. In view of \eqref{fi1}, for \eqref{fi80} it suffices to prove that
\begin{equation}\label{fi81}
||J^{\sigma'}(u_n-u'_n)||_{F^{\sigma_0}}\leq C(\sigma_0,\sigma',||J^{\sigma'}(\phi)||_{H^{\sigma_0}})\cdot ||J^{\sigma'}(\phi-\phi')||_{H^{\sigma_0}},
\end{equation}
for any $n\in\mathbb{Z}_+$. We argue, as before, by induction over $\sigma'$: we decompose $\mathcal{N}(u_n)-\mathcal{N}(u'_n)$ as in \eqref{fi73}, and combine Lemma \ref{Lemmaq1}, Lemma \ref{Lemmaq3}, Lemma \ref{Lemmaq4}, Lemma \ref{Lemmaq6}, and the uniform bound \eqref{fi5}. The proof of \eqref{fi81} is similar to the proof of \eqref{fi20}. This completes the proof of Theorem \ref{Main2}.


\begin{thebibliography}{99}
\bibitem{ChShUh} N.-H. Chang, J. Shatah, and K. Uhlenbeck, Schr\"{o}dinger maps, Comm. Pure Appl. Math, {\bf{53}} (2000), 590--602.
\bibitem{DiWa2} W. Y. Ding and Y. D. Wang, Local Schr\"{o}dinger flow into K\"{a}hler manifolds, Sci. China Ser. A {\bf{44}} (2001), 1446--1464.
\bibitem{IoKe} A. D. Ionescu and C. Kenig, Global well-posedness of the Benjamin--Ono equation in low-regularity spaces, Preprint (2005).
\bibitem{Ka} J. Kato, Existence and uniqueness of the solution to the modified Schr\"{o}dinger  map, Math. Res. Lett., {\bf{12}} (2005), 171--186.
\bibitem{KaKo} J. Kato and H. Koch, Uniqueness of the modified Schr\"{o}dinger map in $H^{3/4+\epsilon}(\mathbb{R}^2)$, Preprint (2005).
\bibitem{KaPo} T. Kato and G. Ponce, Commutator estimates and the Euler and Navier-Stokes equations, Comm. Pure Appl. Math., {\bf{41}} (1988), 891--907.
\bibitem{KeNa} C. E. Kenig and A. Nahmod, The Cauchy problem for the hyperbolic-elliptic Ishimori system and Schr\"odinger maps, Nonlinearity {\bf{18}} (2005), 1987--2009. 
\bibitem{KePoStTo} C. E. Kenig, D. Pollack, G. Staffilani, and T. Toro, The Cauchy problem for Schr\"{o}dinger flows into K\"{a}hler manifolds, Preprint (2005).
\bibitem{KlMa} S. Klainerman and M. Machedon, Space-time estimates for null forms and the local existence theorem, Comm. Pure Appl. Math. {\bf{46}} (1993), 1221---1268.  
\bibitem{KlRo} S. Klainerman and I. Rodnianski, On the global regularity of wave maps in the critical Sobolev norm, Internat. Math. Res. Notices {\bf{13}} (2001), 655--677. 
\bibitem{KlSe} S. Klainerman and S. Selberg, Remark on the optimal regularity for equations of wave maps type, Comm. Partial Differential Equations {\bf{22}} (1997), 901--918. 
\bibitem{Ga} H. McGahagan, An approximation scheme for Schr\"{o}dinger maps, Preprint (2005).
\bibitem{NaStUh} A. Nahmod, A. Stefanov, and K. Uhlenbeck, On Schr\"{o}dinger maps, Comm. Pure Appl. Math., {\bf{56}} (2003), 114--151.
\bibitem{NaStUh2} A. Nahmod, A. Stefanov, and K. Uhlenbeck, Erratum: "On Schr\"odinger maps" [Comm. Pure Appl. Math. {\bf{56}} (2003), 114--151], Comm. Pure Appl. Math. {\bf{57}} (2004), 833--839.
\bibitem{ShSt} J.  Shatah and M. Struwe, The Cauchy problem for wave maps, Int. Math. Res. Notices {\bf{11}} (2002), 555--571.
\bibitem{SuSuBa} P. L. Sulem, C. Sulem, and C. Bardos, On the continuous limit for a system of classical spins, Comm. Math. Phys., {\bf{107}} (1986), 431--454.
\bibitem{Ta1} T. Tao, Global regularity of wave maps. I. Small critical Sobolev norm in high dimension, Internat. Math. Res. Notices {\bf{6}} (2001), 299--328.
\bibitem{Ta2} T.  Tao, Global regularity of wave maps. II. Small energy in two dimensions, Comm. Math. Phys. {\bf{224}} (2001), 443--544.
\bibitem{Tat1} D. Tataru, Local and global results for wave maps. I, Comm. Partial Differential Equations {\bf{23}} (1998), 1781---1793.
\bibitem{Tat2} D. Tataru, On global existence and scattering for the wave maps equation, Amer. J. Math. {\bf{123}} (2001), 37--77.
\bibitem{Tat3} D. Tataru, Rough solutions for the wave maps equation, Amer. J. Math. {\bf{127}} (2005), 293--377.
\end{thebibliography}
\end{document}